\documentclass[12pt]{amsart}

\usepackage{amsmath,amssymb,amscd,amsfonts,eufrak}
\usepackage[mathscr]{eucal}

\newtheorem{thm}{Theorem}[section]
\newtheorem{example}[thm]{Example}
\newtheorem{lem}[thm]{Lemma}
\newtheorem{rem}[thm]{Remark}
\newtheorem{prop}[thm]{Proposition}
\newtheorem{cor}[thm]{Corollary}

\newtheorem{assu-nota}[thm]{Assumption--Notation}
\newtheorem{constr}[thm]{Construction}

\theoremstyle{remark}

\newcommand{\iso}{\cong}
\newcommand{\into}{\hookrightarrow}

\newcommand{\C}{\mathbb C}
\newcommand{\Z}{\mathbb Z}
\newcommand{\Q}{\mathbb Q}

\newcommand{\pp}{\mathbb P}
\DeclareMathOperator{\Aut}{Aut}
\DeclareMathOperator{\Pic}{Pic}
\DeclareMathOperator{\Num}{Num}
\DeclareMathOperator{\Hom}{Hom}
\DeclareMathOperator{\Proj}{Proj}

\newcommand{\OO}{\mathcal O}
\newcommand{\epsi}{\epsilon}

\newcommand{\ga}{\gamma}
\newcommand{\la}{\lambda}
\newcommand{\Ga}{\Gamma}
\newcommand{\De}{\Delta}
\newcommand{\Si}{\Sigma}
\newcommand{\si}{\sigma}
\newcommand{\fie}{\varphi}
\newcommand{\cB}{{\mathcal B}}
\newcommand{\cC}{{\mathcal C}}
\newcommand{\cE}{{\mathcal E}}
\newcommand{\cF}{{\mathcal F}}
\newcommand{\cI}{{\mathcal I}}

\newcommand{\cL}{{\mathcal L}}

\newcommand{\cN}{{\mathcal N}}
\newcommand{\cP}{{\mathcal P}}
\newcommand{\cS}{{\mathcal S}}
\newcommand{\cV}{{\mathcal V}}
\newcommand{\cX}{{\mathcal X}}
\newcommand{\cY}{{\mathcal Y}}
\newcommand{\cZ}{{\mathcal Z}}
\newcommand{\mcE}{\mathscr{E}}
\newcommand{\mcM}{\mathscr{M}}
\newcommand{\tB}{\widetilde{B}}
\newcommand{\tC}{\widetilde{C}}
\newcommand{\tD}{\widetilde{D}}
\newcommand{\tE}{\widetilde{E}}
\newcommand{\tF}{\widetilde{F}}

\newcommand{\tX}{\widetilde{X}}
\newcommand{\tZ}{\widetilde{Z}}
\newcommand{\tpi}{\widetilde{\pi}}
\newcommand{\tq}{\widetilde{q}}
\newcommand{\tcB}{\widetilde{\cB}}
\newcommand{\tcN}{\widetilde{\cN}}
\newcommand{\tcY}{\widetilde{\cY}}
\newcommand{\kS}{\mathfrak{S}}

\newcommand{\inv}{^{-1}}

\numberwithin{equation}{section}

\begin{document}
\title[Surfaces with $p_g=0$ and $K^2=3$]{A new family of surfaces with $p_g=0$ and
$K^2=3$}
\author{Margarida Mendes Lopes \and Rita Pardini}
\date{}

\begin{abstract}  Let  $S$ be a minimal complex surface of general type with $p_g=0$
such that the bicanonical map $\fie$ of $S$ is not birational and let $Z$ be the
bicanonical image. In [M.~Mendes Lopes, R.~Pardini, {\em Enriques surfaces with eight
nodes}, Math. Zeit. {\bf  241} 4 (2002), 673--683] it is shown that either: i)
$Z$ is a rational surface, or ii) $K^2_S=3$, $\fie$ is a degree two morphism and
$Z$ is birational to an Enriques surface. Up to now no example of case ii) was
known. Here an explicit construction of all such  surfaces  is given.  
 Furthermore it is  shown that the corresponding subset of the moduli space of
surfaces of general type is irreducible and uniruled of dimension 6.
\end{abstract}

\maketitle
\section{Introduction} The knowledge of surfaces of general type with $p_g=0$
continues to be scarce  in spite of much progress in surface theory. A minimal
surface of general type  with $p_g=0$ satisfies $1\leq K^2\leq 9$ and examples
for all possible values for $K^2$ are known (see, e.g., \cite{bpv}, Ch. VII, \S 11).

In recent years we have undertaken to study   surfaces with $p_g=0$ by looking at
their bicanonical map, which  is  generically finite onto a surface for  $K^2>1$ (cf. \cite{xiaocan}). When
the bicanonical  map is not birational, this approach works  and it has  allowed
us  to obtain some classification results and also, in some cases, information on
the  moduli space (see for instance
\cite{mp2},
\cite{bica2}, \cite{pianidoppi}).

The first step in describing a class of surfaces with non birational  bicanonical map
is to analyze the bicanonical image.
 In \cite{mp3} the following theorem has been proved:

\begin{thm} \label{8nodes} Let $S$ be a minimal surface of general type such that $p_g(S)=0$,
$K^2_S\ge 3$,  and let
$\fie\colon S\to Z\subset \pp^{K^2_S}$ be the bicanonical map  of $S$. If $\fie$ is
not birational, then either

i) $Z$ is a rational surface,

 or

ii) $K^2_S=3$,
$\fie$ is a morphism of degree 2 and $Z\subset \pp^3$ is an Enriques sextic.

\end{thm} No example of case ii) of the Theorem above appears in the literature.
Indeed,  the known examples of minimal surfaces of general type with
$p_g=0$ and $K_S^2=3$ are the examples of Burniat and Inoue (\cite{bu}, \cite{inoue},
see also \cite{peters}), the examples due independently to J.H. Keum and D. Naie
(\cite{keum}, \cite{naie}) and the recent examples due to F. Catanese
(\cite{quatro}). The degree of the bicanonical map of all these surfaces  is equal to
4 (cf.\cite{mp3},
\cite{mp2}), 
although  the  Keum-Naie examples   are in fact  double covers of
nodal Enriques surfaces and their bicanonical map factorizes through the covering map, as
in case ii) of Theorem \ref{8nodes}.

In this paper not only we show  the existence of surfaces satisfying condition ii),
but  we  give an explicit construction of all such  surfaces
 and we prove  that the corresponding subset of the moduli space of surfaces of
general type is irreducible and uniruled of dimension 6.   The closure of this 
subset contains  the Keum-Naie surfaces (see Proposition \ref{deg4}), whose
fundamental group is known to be isomorphic to $\Z_2^2\times\Z_4$ (cf. \cite{naie}).
Hence the fundamental group of all the surfaces in case ii) of Theorem  \ref{8nodes}
is also isomorphic to $\Z_2^2\times\Z_4$ (see Corollary
\ref{fgroup}).

\vskip0.2truecm Our  description   of surfaces satisfying condition ii) is based on  
 a very detailed study of the normalization of their  bicanonical
images. These are  
polarized Enriques surfaces  of degree 6 with 7 nodes, satisfying some
additional conditions (see Proposition
\ref{gentoen}  and the setting of \S \ref{secex}). The analysis and construction of
these Enriques surfaces  form the bulk of this paper. The main tools we use are the
classification of linear systems on an Enriques surface, the analysis of the 
configuration of singular fibres of certain elliptic pencils, the code associated to
the nodes of the surface and the corresponding Galois cover (cf. \cite{nodes}).

\vskip0.2truecm The plan of the paper is as follows: Section~\ref{seconda} explains
the relation between the surfaces satisfying condition ii) of Theorem \ref{8nodes} and a
certain class of polarized Enriques surfaces $(\Si, B)$ with 7 nodes; in
Section~\ref{secex} some properties of these Enriques surfaces are established and 
 some examples are described; in  Section~\ref{seccode} we make a very detailed study of the
singular  fibres of the elliptic pencils of $\Si$ and we determine  the code associated to
the nodes of $\Si$; in Section~\ref{secconstr} we describe a construction yielding pairs
$(\Si, B)$  and prove that all such pairs are obtained in that way; in Section~\ref{secmod}
we introduce and
 study a quasi-projective  variety parametrizing the isomorphism classes of pairs
$(\Si,B)$ and finally in Section~\ref{newfamily} we apply the previous results to
describe the family of surfaces satisfying condition ii) of Theorem \ref{8nodes}.

\vskip0.3truecm  \noindent {\bf Notation and conventions:} 
We work over the complex numbers.
A node of a surface is an ordinary double point, namely a singularity  analytically
isomorphic to $x^2+y^2+z^2=0$. The exceptional divisor of the minimal resolution of a node
is a curve  $C\simeq \pp^1$ such that $C^2=-2$. A curve with these properties is called 
$-2-$curve, or  nodal curve.

 We say that a projective surface $\Si$ with canonical singularities
is minimal, of general type, Enriques
\dots if the minimal resolution of $\Si$ is minimal, of general type, Enriques \dots.
Our standard reference for Enriques surfaces is \cite{codo}, and we use freely its 
terminology.

Given an automorphism $\si$ of a variety $X$, we say that a map $f\colon X\to Y$ is composed with
$\si$ if $f\circ\si=f$.   If $G$ is a finite group, a $G-$cover  is a finite map $f\colon X\to Y$ of
normal varieties together with a faithful $G-$action on $X$ such that $f$ is isomorphic to
the quotient map
$X\to X/G$. If $G=\Z_2$, then  we say that $f$ is a double cover. Contrary to what is often
done (cf. for instance 
\cite{ritaabel}), we do not require that
$f$ be flat.

We denote linear equivalence by $\equiv$ and numerical equivalence by~$\sim_{num}$. The group of line
bundles modulo numerical equivalence on a variety $Y$ is denoted by $\Num(Y)$.
\medskip

\noindent{\bf Acknowledgements.} 
 The present collaboration takes place in the framework of the european contract
EAGER, no. HPRN-CT-2000-00099. It was partially supported by the 1999 Italian P.I.N.
``Geometria sulle Variet\`a Algebriche'' and by the ``Financiamento Plurianual'' of
CMAF. The first author is a member of CMAF and of the Departamento de Matem\'atica da
Faculdade de Ci\^encias da Universidade de Lisboa and the second author is a member
of GNSAGA of CNR.  

We are indebted to Igor Dolgachev, Barbara Fantechi and Marco Manetti for useful discussions whilst writing this paper.

\section{Surfaces with $p_g=0$ and $K^2=3$ and Enriques surfaces}\label{seconda} As
explained in the introduction, the bulk of this paper is   a very detailed  study of
a class of polarized Enriques surfaces with 7 nodes. 
 In this section we explain the relation  between such Enriques surfaces
 and  a  class of minimal surfaces of general type with non birational bicanonical
map.  Let
$S$ be  a minimal surface of general type with
$p_g(S)=0$ and
$K^2_S=3$. We denote by $\fie\colon S\to\pp^3$ the bicanonical map of $S$ and we
assume that $S$ has an involution $\si$ such that:

a) $\fie$ is composed with $\si$;

b) the quotient surface
$T:=S/\si$ is birational to an Enriques surface.

We denote by $X$ the canonical model of $S$. Abusing notation, we denote by the same
letter the involution induced  by
$\si$  on $X$.

\begin{prop}\label{gentoen}In the above setting:
\begin{enumerate}
\item the quotient surface $\Si:=X/\si$ is an Enriques surface with 7 nodes;
\item the quotient map $\pi\colon X\to\Si$ is branched on the nodes of
 $\Si$ and on a divisor $B$  with negligible singularities,  contained in the smooth
part of $\Si$;
\item $B$ is ample and $B^2=6$;
\item the bicanonical system $|2K_X|=\pi^*|B|$ is base point free;
\item the degree of the bicanonical map of $S$ (and $X$) is either 2 or 4.
\end{enumerate}
\end{prop}
\begin{proof} Since the bicanonical map $\fie$ of $S$ factorizes through $\si$, by
\cite[Prop. 2.1]{mp2} the isolated fixed points of $\si$ are 7. The quotient surface
$T:=S/\si$ has 7 nodes, which are the images of the isolated fixed points of $\si$.
The quotient map $S\to T$ is branched on the nodes and on a smooth  divisor $B_0$
contained in the smooth part of $T$.  By Lemma 7 of
\cite{xiao2} and the Remark following it, there exists a birational morphism
$r\colon T\to \Si'$, where $\Si'$ is an Enriques surface with 7 nodes, such that the
exceptional curves of $r$ are contained in the smooth part of $T$ and the divisor
$B:=r(B_0)$ has negligible singularities. Let $S\to X'\overset{\pi'}{\to}  \Si'$ be
the Stein factorization of the induced map $S\to\Si'$. The map $\pi'\colon X'\to\Si'$
is a double cover branched on the nodes of $\Si$ and on the divisor $B$. The
singularities of $X'$ occur above the singularities of $B$, hence they are canonical
and there is a birational morphism from $X'$ to the canonical model
$X$  of $S$. More precisely, there is a commutative diagram:
\[
\begin{CD} X'@>>>X\\ @V{\pi'}VV @VV{\pi}V\\
\Si'@>>>\Si
\end{CD}
\] where the horizontal arrows represent birational morphisms and $\pi$, $\pi'$ are
the quotient maps for the involutions induced by $\si$ on $X$ and $X'$.

 By adjunction we have
$2K_{X'}={\pi'}^*(2K_{\Si'}+B)={\pi'}^*B$, hence
$B$ is nef and $B^2=6$.
 Since the bicanonical map of
$X'$ factorizes through
$\pi'$, we actually have
$|2K_{X'}|={\pi'}^*|B|$. The same argument as in the proof of \cite[Thm. 5.1]{mp3}
shows that the system $|B|$ is free. Thus  $(\Si',B)$ is a pair as in the setting of \S
\ref{secex} and   we can apply  Corollary
\ref{Bample}, which is proven in \S \ref{secex}, showing that  $B$ is ample. It follows that  the
horizontal maps in the above diagram are isomorphisms and we can identify $X$ with
$X'$ and $\Si$ with $\Si'$.  The system  $|B|$, being free,  is not hyperelliptic,
therefore   by \cite[Prop. 5.2.1]{cossec} either it is birational or it has degree 2.
Hence the degree of the bicanonical map of $X$ and  $S$ has  degree either 2 or 4.
This completes the proof.
\end{proof}

The previous Proposition has a converse:
\begin{prop} \label{entogen}Let $\Si$ be an Enriques surface with 7 nodes and let
$\pi\colon X\to \Si$ be a
 double cover branched on the nodes of $\Si$ and on a divisor $B$ such that:

a) $B$ is ample and $B^2=6$;

b) $B$ is contained in the smooth part of $\Si$ and it has negligible singularities.

Then $X$ is the canonical model of a minimal surface $S$ of general type with
$p_g(S)=0$ and
$K^2_S=3$ and the bicanonical map of $S$ factorizes through the map $S\to \Si$.
\newline Condition a) can be replaced by:

a') $|B|$ is free and $B^2=6$.
\end{prop}
\begin{proof} Assume  that conditions a) and b) are satisfied.\newline  The singular
points of $X$ lie above the singularities of $B$. Since $B$ has negligible
singularities, the singularities of $X$ are canonical and one has $2K_X=\pi^*B$. 
Since $B$ is ample, $K_X$ is also ample and $X$ is the canonical model of a surface
$S$  of general type. One has $K^2_S=K_X^2=\frac{1}{2}B^2=3$. To compute the
birational invariants $\chi(S)$ and $p_g(S)$, we consider the minimal resolution of  
singularities  $\eta\colon Y\to\Si$ and  the flat  double cover
$\tpi\colon
\tX\to Y$ obtained from $\pi$ by taking base change with $\eta$. The surface $\tX$ 
has  canonical singularities and it is birational to $X$ and $S$. The branch locus of
$\tpi$ consists of the inverse image $\tB$ of $B$ and of the  $-2-$curves
$N_1,\dots, N_7$ that are exceptional  for $\eta$. So $\tpi$ is given by a relation
$2L\equiv
\tB+N_1+\dots+N_7$, where $L$ is a line bundle on $Y$. A standard computation gives
 $\chi(S)=\chi(\tX)=\chi(Y)+\chi(L\inv)=1$.  By Kawamata--Viehweg vanishing one has
$h^i(K_Y+L)=0$ for $i>0$, hence $p_g(S)=p_g(\tX)=h^0(K_Y+L)=\chi(K_Y+L)=0$.

Since the singularities of the branch divisor $B$ of $\pi$ are negligible, it is well
known that the smooth minimal  model $S$ of $X$ can be obtained by repeatedly blowing
up $\Si$ at the singular points of $B$ and taking base change and normalization.
Hence the involution of $S$ induced by $\si$ has the same number  of isolated base
points as $\si$, that is 7,  and the bicanonical map factorizes through it by
\cite[Prop. 2.1]{mp2}.

If a') holds, then $B$ is ample by Corollary \ref{Bample}, and so  condition a')
implies condition a).
\end{proof}

\section{Enriques surfaces with 7 nodes:  examples}\label{secex}
Recall that for Enriques surfaces we adopt the notation and the terminology of \cite{codo}. Our
notation for the singular fibres of an elliptic pencil is the same as in \cite{bpv}, Ch. V, \S 7. 

In this section we  consider  the following situation:

\noindent{\bf Set--up:} $\Si$ is a nodal Enriques surface with 7 nodes,
$|B|$ is a base point free linear system of $\Si$ such that $B^2=6$. We denote by
$\eta\colon Y\to\Si\/$\  the minimal desingularization and by
$N_1,
\dots, N_7$ the disjoint nodal curves contracted by $\eta$, and we set $\tB:=\eta^*B$.
Furthermore, we assume that there exists
$L\in
\Pic(Y)$ such that $\tB+N_1+\dots+N_7=2L$.

\begin{rem}{\em The choice of this set--up is suggested by the results of the
previous section. Indeed, the condition that the class
$\tB+N_1+\dots +N_7$ be divisible by 2 in
$\Pic(Y)$  means that, given a curve
$\tB\in |\tB|$, there exists a double cover
$\tpi\colon \tX\to Y$ branched on the union of $\tB$ and $N_1,\dots, N_7$. If $\tB$ is
disjoint from
$N_1,\dots ,N_7$ (recall that $\tB N_i=0$) and it has at most negligible singularities,
then the surface
$\tX$ has canonical singularities, occurring above the  singularities of $\tB$. For
$i=1,\dots,7$, one has
$\tpi^*N_i=2C_i$, where $C_i$ is a $-1-$curve contained in the smooth part of $\tX$.
If we denote by 
$X$ the surface obtained by contracting the $C_i$, then $\tpi$ induces a double 
cover $\pi\colon X\to\Si$ branched over the image $B$ of $\tB$ and over the nodes of
$\Si$. By  Proposition
\ref{entogen}, $X$ is the canonical model of a minimal surface of general type with
$p_g=0$ and
$K^2=3$ and the bicanonical map of $X$ factorizes through $\pi$.

Notice that, since $Y$ is an Enriques surface, the line bundle $L+K_Y$ also satisfies
the relation $2(L+K_Y)\equiv \tB+N_1+\dots +N_7$, so that a pair $(\Si, B)$ as in the
set--up determines two non-isomorphic double covers of $\Si$ with the same branch
locus.}
\end{rem}

In order to describe some examples, we need to prove first some general facts.

\begin{prop}\label{3pencils} There exist three elliptic half-pencils
$\tE_1,\tE_2,\tE_3$ on $Y$ such that:
\begin{enumerate}
\item $\tE_i\tE_j=1$ for $i\ne j$;
\item $|\tB|=|\tE_1+\tE_2+\tE_3|$
\end{enumerate}
\end{prop}
\begin{proof} Notice that the system $|\tB|$, being base point free, is not
hyperelliptic.  Hence, if
$|\tB|$ is not as stated, then by Proposition 5.2.1 and Theorem 5.3.6 of
\cite{cossec} there are the following possibilities:
\begin{itemize}
\item[1)] $|\tB|=|2\tE_0+\tE_1+\theta_2|$, $\tE_0\tE_1=\tE_0\theta_2=1$,
$\tE_1\theta_2=0$ 
\item[2)] $|\tB|=|3\tE_0+2\theta_0+\theta_1|$,
$\tE_0\theta_0=\theta_0\theta_1=1$,$\tE_0\theta_1=0$
\end{itemize} where   $\tE_m$ are  elliptic half-pencils and  $\theta_m$  are
nodal curves, $m=0,1,2$. Consider the  nodal curves $N_1,\dots, N_7$
and recall that
$\tB+N_1+\dots+ N_7$ is divisible by 2 in $\Pic(Y)$.

In case 1), suppose first that $\theta_2$ is not one of the curves $N_i$.  Since
$N_i\tB=0$, necessarily the curves $\theta_2, N_1, \dots, N_7$ are disjoint and so by Lemma 4.2
of \cite{mp3}  the divisor $\theta_2+ N_1+ \dots+ N_7$ is divisible by 2 in
$\Pic(Y)$. Hence also
$\tE_1$ is divisible by 2 in  $\Pic(Y)$, a contradiction. So $\theta_2$ is one of the
curves $N_i$, say $\theta_2=N_7$. Then there exists $L\in \Pic(Y)$ such that
$D=\tE_1+N_1+\dots+N_6= 2L$ and we get $D^2=-12$ and $L^2=-3$, contradicting the fact
that the intersection form on an Enriques surface is even.

Consider now case 2). As in case 1), if neither of the curves
$\theta_0,\theta_1$ is one of the nodal curves $N_1, \dots, N_7$, we conclude that
the half-pencil $\tE_0$ is divisible by 2 in $\Pic(Y)$, a contradiction.  If
$\theta_1$ is one of the curves
$N_i$, we also arrive at a contradiction as in case 1). Finally suppose that
$\theta_0$ is one of the curves $N_i$, say $N_7$. Then $\tE_0(N_1+ \dots+ N_7)=1$.
Since
$\tE_0\tB=2$,
$\tB+N_1+\dots+ N_7$ is not divisible by 2 in $\Pic(Y)$, a contradiction.

So $|\tB|$ is as stated.
\end{proof}

\begin{cor}\label{Bample} The divisor $B$ is ample on $\Si$.
\end{cor}
\begin{proof} Denote by ${\mathcal R}_{\tB}$ the set of irreducible curves $C$ of
$Y$  with $\tB C=0$.  By \cite[Cor. 4.1.1]{codo}, ${\mathcal R}_{\tB}$ is a  finite
set and the corresponding classes are independent in $\Num(Y)$.    By Proposition
\ref{3pencils}, we can write
$\tB=\tE_1+\tE_2+\tE_3$, where the
$\tE_i$ are elliptic half--pencils.  If we denote by $V$ the subspace of $\Num
(Y)\otimes \Q$ spanned by the classes of 
$\tE_1, \tE_2,
\tE_3$, then the classes of the curves of
${\mathcal R}_{\tB}$ belong to $V^{\perp}$, which has dimension 7. So we have
 ${\mathcal R}_{\tB}=\{N_1,\dots ,N_7\}$ and 
$B=\eta_*\tB$ is ample on $\Si$.
\end{proof}

The following is a partial converse to Proposition \ref{3pencils}:
\begin{lem}\label{even} Let $Y$ be a smooth Enriques surface containing 7 disjoint
nodal curves
$N_1,\dots ,N_7$.  Assume that
$\tE_1,
\tE_2, \tE_3$ are elliptic half pencils on $Y$ such that $\tE_i\tE_j=1$ if
$i\ne j$ and $\tE_iN_j=0$ for every $i,j$. If we  set $\tB:=\tE_1+\tE_2+\tE_3$, then
$\tB+N_1+\dots +N_7$ is divisible by 2 in
$\Num(Y)$.
\end{lem}
\begin{proof} Recall that $\Num(Y)$ is an even unimodular lattice of rank 10.
\newline Let $M$ be the sublattice of $\Num(Y)$ spanned by the classes of
$\tE_1-\tE_2$ and
$\tE_1-\tE_3$. The discriminant of $M$ is equal to 3, hence $M$ is primitive.
 Let $M'$ be the sublattice spanned by the classes of the
$N_i$ and by the class of $\tB$. The primitive closure of $M'$ is
$M^{\perp}$ and the code $W$ associated to the set of classes $\tB, N_1,\dots ,N_7$
is naturally isomorphic to the quotient group
$M^{\perp}/M'$.  Computing discriminants one gets:  $$2^8\cdot 3=
disc(M')=disc(M^{\perp})2^{2 dim W}=disc(M)\cdot 2^{2dim W}=3\cdot 2^{2dim W},$$
namely $\dim W=4$. Using the fact that the intersection form on $Y$ is even, it is
easy to check that the elements of
$W$ have weight divisible by 4.   Since $W$ has length 8, this implies that
$W$ is the extended Hamming code (see, e.g., \cite{lint}), and in
particular it contains the vector of weight 8, i.e.,
$\tB+N_1+\dots +N_7$ is divisible by 2 in $\Num(Y)$.
\end{proof}
 For $i=1,2,3$, we denote by $\tE_i'$ the unique effective divisor in $|\tE_i+K_Y|$
and we write 
$|\tF_i|=|2\tE_i|=|2\tE_i'|$. Thus   $|\tF_i|$ is an elliptic pencil with  double
fibres
$2\tE_i$ and $2\tE_i'$.  The classes $\tE_i$ are nef, hence $\tB N_j=0$ implies $\tE_i
N_j=0$ for every $i,j$ and therefore for
$i=1,2,3$
$|\tF_i|$ induces an elliptic pencil $|F_i|=|2E_i|=|2E_i'|$ on $\Si$, where
$E_i=\eta_*\tE_i$, $E_i'=\eta_*\tE_i'$ and
$B=E_1+E_2+E_3$.

\begin{example}\label{ex1}{\em  This example appears in \cite{naie} and in an
unpublished paper by J. Keum (\cite{keum}). One  considers an Enriques surface
$\overline{\Si}$ with 8 nodes as in Example 1 of
\cite{mp3}. The surface $\overline{\Si}$ has two isotrivial elliptic pencils $|F_1|$
and
$|F_2|$ with
$F_1 F_2=4$. The system $|F_1+F_2|$ gives a degree 2 morphism onto a Del Pezzo
quartic in
$\pp^4$ such that the nodes of $\overline{\Si}$ are mapped to smooth points (cf.
\cite[\S 2]{naie}). We  take
$\Si$ to be the surface obtained 
 by resolving one of the nodes $\overline{\Si}$ and we denote by $C$ the corresponding
nodal curve.  We denote by the same letter the pull--backs  of 
$|F_1|$, $|F_2|$ on $\Si$  and we set
$B:=F_1+F_2-C$.  By the above discussion, the system $|B|$ is free  and  it gives a 
degree 2 map onto a Del Pezzo cubic in $\pp^3$.    By
\cite[Lemma 4.2]{mp3}, the class of $N_1+\dots +N_7+C$ is divisible by 2 in
$\Pic(Y)$, hence the class of
$\tB+N_1+\dots+N_7$ is also divisible by 2.

Let $2E_i$ be a double fibre of $F_i$, $i=1,2$. By Riemann--Roch there exists an effective divisor 
$E_3\equiv E_1+E_2-C$.  We have
$E_3^2=0$ and
$E_1E_3=E_2E_3=1$. We claim that $E_3$ is an elliptic half-pencil, so that $B\equiv E_1+E_2+E_3$ as
predicted by Proposition \ref{3pencils}. 

 We now work on the non singular surface $Y$ and, as usual,  we denote by
$\tD$ the pull back on $Y$ of a divisor $\tD$ of $\Si$.  Let
$G_i\in |\tF_i|$,
$i=1,2$, be the fibre containing 
$C$. By the description of $\overline{\Si}$ given in \cite{mp3}, $G_1$ and $G_2$ are fibres of type
$I_0^*$ and the divisor $G_3:=G_1+G_2-2C$ is an elliptic configuration of type $I_2^*$. It follows
that
$|G_3|$ is an elliptic pencil and that  $2\tE_3$ is a double fibre of $|G_3|$.}
\end{example}

\begin{example}\label{ex2} {\em Let ${\cC}\subset \pp^3$ be the Cayley cubic, defined
by
$x_1x_2x_3+x_0x_2x_3+x_0x_1x_3+x_0x_1x_2=0$. The singularities  of $\cC$ are 4 nodes,
that occur at the coordinate points and form an even set. The 6 lines joining the
nodes are of course contained in
$\cC$. We label these lines by $e_1, e_1', e_2, e_2' , e_3, e_3'$ in such a way that
for $i=1,2,3$
$e_i, e_i'$ is a pair of skew lines and $e_1, e_2, e_3'$ are coplanar. The surface
$\cC$ contains 3 more lines
$l_1, l_2, l_3$, contained in the plane
$x_0+x_1+x_2+x_3=0$. An elementary geometric argument shows that, up to a permutation
of the indices, we may assume that
 the line $l_i$ meets $e_i,e'_i$ and it  does not meet $e_j,e'_j$ for $i\ne j$. For
$i=1,2,3$ we denote by $|f_i|$ the moving part  of the linear system cut out on $\cC$
by the planes containing the line $l_i$. The general $f_i$ is a smooth conic, hence
one has
$K_{\cC}f_i=-2$, $f_i^2=0$ and $f_if_j=2$ for $i\ne j$. The singular fibres of
$|f_i|$ are
$2e_i$,
$2e_i'$ and
$l_j+l_k$, where $i,j,k$ is a permutation of $1,2,3$.

  Consider a curve $D\in |\OO_{\cC}(2)|$  such  that
$D$ is contained in the smooth part of $\cC$ and it  has at most simple singularities.
Let $\Si$ be the double cover of $\cC$ branched on $D$ and on   the four nodes of
$\cC$. The surface $\Si$ has canonical singularities, occurring over the singular
points of $D$. Standard computations (cf. the proof of Proposition
\ref{entogen}) show that
$\Si$ is an Enriques surface and that the  pull back of the system of hyperplanes  of
$\pp^3$ is a complete system $|B'|$ on $\Si$ with
${B'}^2=6$.  For $i=1,2,3$, we consider on $\Si$ the system $|F_i|$ obtained by
pulling back $|f_i|$. The system $|F_i|$ is an elliptic pencil, with double fibres
$2E_i$ and $2E_i'$, where $E_i$, $E_i'$ are the pull backs of $e_i$, and  $e_i'$,
respectively. For
$i\ne j$ one has $E_iE_j=1$. Furthermore, $B'\equiv E_1+E_2+E_3'$.
\bigskip

\begin{center}
\setlength{\unitlength}{0.5truecm}
\begin{picture}(16,15)\label{I21}
\put(8,2){\circle*{0.25}}
\put(8,13){\circle*{0.25}}
\put(8,5){\circle*{0.25}}
\put(4,7){\circle*{0.25}}
\put(12,7){\circle*{0.25}}
\put(4,11){\circle*{0.25}}
\put(12,11){\circle*{0.25}}
\put(1,13){\circle*{0.25}}
\put(15,13){\circle*{0.25}}
\thicklines
\put(8,2){\line(0,1){3}}
\put(8,5){\line(2,1){4}}
\put(8,5){\line(-2,1){4}}
\put(4,7){\line(0,1){4}}
\put(12,7){\line(0,1){4}}
\put(4,11){\line(-3,2){3}}
\put(4,11){\line(2,1){4}}
\put(12,11){\line(-2,1){4}}
\put(12,11){\line(3,2){3}}
\put(7.5,5.8){$A_1$}
\put(12.5,10.5){$A_2$}
\put(2.5,10.5){$A_3$}
\put(6.5,1.5){$N_4$}
\put(2.5,6){$N_2$}
\put(12.5,6){$N_3$}
\put(7,13.5){$N_1$}
\put(0.8,13.5){$N_6$}
\put(15,13.5){$N_5$}
\end{picture}

Figure \ref{I21}
\end{center}
\bigskip

We consider now a special case of the above construction: we take $D$ to be the union
of the section
$H_0$ of $\cC$ with the plane $x_0+x_1+x_2+x_3=0$ and of the section $H$ with a
general hyperplane tangent to
$\cC$. So $H$ has an ordinary double point at the tangency point and is smooth
elsewhere and $H$ and
$H_0$ intersect transversely at 3 points. The surface $\Si$ thus obtained has 7
nodes, occurring above the singularities of $H_0+H$. As usual we denote by
$\eta\colon Y\to
\Si\,$ the minimal resolution and by $N_1,\dots ,N_7$ the exceptional curves of
$\eta$. The strict transform
 on
$Y$  of the line
$l_i$ is a nodal curve $A_i$. It is not difficult to check that one can relabel the curves
$N_1,\dots ,N_7$ in such a way that $N_7$ corresponds to the  singularity of $\Si$
above the double point of
$H$ and the incidence relations of the set of curves $A_1,A_2, A_3,N_1,\dots, N_6$
are as shown in the dual graph in Figure \ref{I21}}

{\em As usual we denote by  $\tE_i$, $\tE_i'$, $\tF_i$, $\tB'$ the pull backs on $Y$
of $E_i$, $E_i'$,
$F_i$, $B'$. The singular fibre  $G_i$  of
$|\tF_i|$  corresponding to the fibre $l_j+l_k$ of $|f_i|$ is  of type
$I_2^*$. More precisely, we  have $G_1=N_2+N_3+N_5+N_6+2(A_3+A_2+N_1)$,
$G_2=N_1+N_3+N_4+N_6+2(A_3+A_1+N_2)$,
$G_3=N_1+N_2+N_4+N_5+2(A_1+A_2+N_3)$. By Lemma \ref{even}, one of the classes
$B'+N_1+\dots+N_7$ and
$B'+K_Y+N_1+\dots +N_7$ is divisible by
$2$ in $\Num (Y)$. We will show later (Corollary \ref{corex2}) that the second case
actually occurs. Hence we set
$B:=B'+K_{\Si}$.}
\end{example}

Examples \ref{ex1} and \ref{ex2} share the common feature that either the system
$|B|$ or the system
$|B+K_{\Si}|$ is not birational. The next example shows that this does not happen in
general.

\begin{example}\label{ex3} {\em By deforming Example 1 we show the existence of  a
pair  $(\Si, B)$  such that both $|B|$ and
$|B+K_{\Si}|$ are birational .

We start with  a pair
$(\Si_0, B_0)$ as in  Example \ref{ex1}. The Kuranishi family  $p\colon
\cY\to U$ of the minimal resolution
$Y_0$ of $\Si_0$ is smooth of dimension 10 by \cite[Thm. VIII.19.3]{bpv}. We may
assume that $U$ is contractible and that the family
$\cY$ is differentiably trivial. Hence for every fibre $Y_t:=p\inv(t)$ the inclusion
$Y_t\to \cY$ induces an isomorphism
$H^2(Y_t,\Z)\stackrel{\sim}{\longrightarrow}H^2(\cY,\Z)$.  The Leray spectral sequence
gives
$h^1(\cY,\OO_{\cY})=h^2(\cY,\OO_{\cY})=0$, hence by the exponential sequence every
integral cohomology class of
$\cY$ comes from a unique holomorphic line  bundle. In particular, there exist line
bundles
$\widetilde{\cE_1}, \widetilde{\cE_2}, \widetilde{\cE_3}, \cN_1, \dots,
\cN_7,\cC\/$ that restrict on the central fibre
$Y_0$ to
$\tE_1,\tE_2,\tE_3$, $N_1,\dots, N_7, C$, respectively. We set
$\widetilde{\cB}:=\widetilde{\cE_1}+\widetilde{\cE_2}+\widetilde{\cE_3}$ and for
$t\in U$ we  denote by
$\tB_t$, $\tE_{i,t}$, $N_{i,t}$,
$C_t$ the restrictions to $Y_t$ of the above bundles. Obviously, the class of $\tB_t+N_{1,t}+\dots
+N_{7,t}$ is divisible by 2 in
$H^2(Y_t,\Z)$ for every $t\in U$.

By  \cite[Theorem 3.7]{bw}, the subset $U_1$ of $U$ where the classes
$N_{i,t}$ are effective and irreducible is smooth of dimension 3, while the  subset
$U_2$ of $U_1$ where also $C_t$ is effective is smooth of dimension $2$. Since
$C_t\equiv \tE_{1,t}+\tE_{2,t}-\tE_{3,t}$, by
\cite[Thm. 4.7.2]{codo}  the system
$\tB_t$ is birational for $t\in U_1\setminus U_2$.  On the other hand, by
semicontinuity we may assume that
$\tB_t+K_{Y_t}$ is birational for every $t\in U$, since it is birational on the
central fibre $Y_0$. So the required example can be obtained by taking $Y_t$ with
$t\in U_1\setminus U_2$ and by blowing down the nodal curves $N_{1,t},\dots, N_{7,t}$.}
\end{example}

\section{Enriques surfaces with 7 nodes: codes and singular fibres}\label{seccode}

We keep the set--up and the notation of the previous section.
\newline Here we make a detailed study of the code associated to the nodal curves
$N_1,\dots, N_7$ and of the singular fibres of the pencils $|\tF_i|$. These results
are needed in the following section, where  we give a construction of all the pairs
$(\Si,B)$ as in the set--up of \S  \ref{secex}.

We denote by  $V$  and $V_{num}$, respectively, the code and the numerical code
associated to
$N_1,\dots, N_7$ (cf. \cite[\S 2]{mp3}). Namely, $V$ is the kernel of the map
$\Z_2^7\to \Pic(Y)/2\Pic(Y)$ that maps $(x_1,\dots,x_7)$ to the class of
$x_1N_1+\dots+x_7N_7$. The code   $V_{num}$  is defined in analogous way, replacing
$\Pic(Y)$ by $\Num(Y)$. Clearly,  $V$ is a subcode of $V_{num}$ of codimension
$\le 1$.  We say that a divisor  $D$ is even if it is divisible by 2 in $\Pic(Y)$. In
particular, if $D=\sum x_iN_i$ then $D$ is even if and only if 
$(x_1,\dots, x_7)\in V$ (we denote by the same letter the integer $x_i$ and its
class in
$\Z_2$). 
\begin{lem}\label{Vnum}
$\dim V_{num}=3$.
\end{lem}
\begin{proof} Since the determinant of the matrix $(\tE_i\tE_j)_{i,j=1\dots 3}$ is
equal to 2, the classes $\tE_1, \tE_2,\tE_3$ span a primitive sublattice $L$  of rank
3 of
$\Num(Y)$. If $L'$ is the sublattice spanned by the classes $N_1,\ldots N_7$, then
$V_{num}$ is isomorphic to the quotient group $L^{\perp}/L'$. So we have
$$2^7=disc(L')=2^{2\dim V_{num}}disc(L^{\perp})=2^{2\dim V_{num}}disc(L),$$ namely
$\dim V_{num}=3$.

\end{proof}

\begin{lem}\label{Bsepara} The linear system
$|B|$ separates the nodes of $\Si$.
\end{lem}
\begin{proof} This follows by  \cite[Lemma 4.6.3]{codo}.
\end{proof}
\begin{lem}\label{collinear} Denote by $P_1,\dots, P_7$ the image points
  of
$N_1,\dots, N_7$ via the system $|\tB|$. If
$N_1+N_2+N_3+N_4$ is an even divisor, then $P_5,P_6,P_7$ are collinear.
\end{lem}
\begin{proof} Notice first of all that the points $P_1,\dots, P_7$ are distinct by
Lemma
\ref{Bsepara}.

 Since $\tB+N_1+\dots +N_7$ is even by assumption, there exists $M\in
\Pic(Y)$  such that
$2M\equiv \tB-N_5-N_6-N_7$. Set $M':=M+K_Y$. Since $M^2={M'}^2=0$, there exist
effective divisors $D\in |M|$ and $D'\in |M'|$. So the linear system
$|\tB-N_5-N_6-N_7|$ contains two distinct divisors $2D$ and $2D'$, hence it has
positive dimension. This means that
$P_5, P_6,P_7$ lie on a line.
\end{proof}

\begin{prop}\label{dimV} The code $V$ has dimension 2.
\end{prop}
\begin{proof} By Lemma \ref{Vnum}, to show that $\dim V=2$ it is enough to show that
$V\subsetneq V_{num}$. So assume by contradiction that $V=V_{num}$. Since
$\dim V=3$ and all the elements of
$V$ have weight 4,  $V$ is isomorphic to the Hamming code (see, e.g., \cite{lint}).  By the definition of the Hamming code,
the  set of indices $\{1,\dots, 7\}$ is in one-to-one correspondence with the nonzero vectors of
$\Z_2^3$. The vectors corresponding to distinct indices $i_1$, $i_2$, $i_3$ span a
plane  of
$\Z_2^3$ if and only if there is  $v=(x_1,\dots, x_7)\in V\setminus\{0\}$ such that
$x_{i_1}=x_{i_2}=x_{i_3}=0$. By Lemma \ref{collinear}, this happens if and only if
the points
$P_{i_1},P_{i_2},P_{i_3}$ lie on a line in $\pp^3$. Hence the points
$P_1,\dots, P_7$ form a configuration isomorphic to the finite plane $\pp^2(\Z_2)$.
Since the line through two of the
$P_i$ contains a third point of the set, it is easy to check that $P_1,\dots, P_7$ lie
in a plane. On the other hand, it is well known that the plane $\pp^2(\Z_2)$ cannot be
embedded in $\pp^2(\C)$. So we have a contradiction and the proof is complete.
\end{proof} We are now able to complete the description of Example \ref{ex2}:
\begin{cor}\label{corex2} Let $\Si$ be the surface of Example \ref{ex2}. Then
$\tB'+K_Y+N_1+\dots+ N_7$ is divisible by 2 in $\Pic(Y)$.
\end{cor}
\begin{proof} In  the notation of Example \ref{ex2} we have
$G_1=2(A_2+A_3+N_1)+N_2+ N_3+N_5+N_6\equiv 2\tE_1$, hence $N_2+N_3+N_5+N_6$ is an
even divisor. The same argument shows that the divisors $N_1+N_2+N_4+N_5$ and
$N_1+N_3+N_4+N_6$ are also even. By Proposition \ref{dimV} these are the only non
zero elements of $V$.

By Lemma \ref{even} we know that one of the classes $\tB'+N_1+\dots+N_7$ and
$\tB'+K_Y+N_1+\dots+N_7$ is even. Assume by contradiction that
$\tB'+N_1+\dots+N_7$ is even. Pulling back to $Y$ the section of $\cC$ with $H_0$ we
get $\tB'\equiv 2(A_1+A_2+A_3+N_1+N_2+N_3)+N_4+N_5+N_6$. Hence it follows that
$N_1+N_2+N_3+N_7$ is also an even divisor, a contradiction.
\end{proof}

The next result describes the possible  configurations of  singular fibres of the
pencils
$|\tF_i|$ and relates them to the properties of the systems $|B|$ and
$|B+K_{\Si}|$. 

\begin{thm}\label{B+fibres} The possible  configurations of fibres with singular
support of the pencils $|\tF_i|$ are the following:
\begin{itemize}
\item[1)] up to a permutation of the indices, the pencils $|\tF_1|$ and $|\tF_2|$ are
isotrivial with 2 fibres of type $I_0^*$, while $|\tF_3|$ has a fibre of type $I_2^*$
and two fibres of type
$I_2$ or $_2I_2$.
\newline In this case the system $|B|$ has degree 2 and the system $|B+K_{\Si}|$ is
birational.
\item[2)] each of the pencils $|\tF_i|$ has a fibre of type $I_2^*$ and two fibres of
type
$I_2$ or $_2I_2$. The dual graph of the set of nodal curves that form  the
$I_2^*$ fibres is the same as in Figure \ref{I21}.
\newline In this case $|B|$ is birational and $|B+K_{\Si}|$ has degree 2.

\item[3)] each of the pencils $|\tF_i|$ has a fibre of type $I_0^*$ and three fibres
of type
$I_2$ or $_2I_2$. 
\newline In this case the systems $|B|$ and $|B+K_{\Si}|$ are both birational.

\end{itemize}
\end{thm}
\begin{rem}{\em The proof of Theorem  \ref{B+fibres} below actually  shows more,
namely  that  case 1) of Theorem
\ref{B+fibres} corresponds exactly to Example \ref{ex1} (cf. Lemma \ref{extranode})
and that case 2) corresponds exactly to Example \ref{ex2}.}
\end{rem} The proof of Theorem \ref{B+fibres} is somewhat involved and requires some
auxiliary lemmas.
\begin{lem}\label{scemo}  Assume that for
$i\ne j$, the pencils $|\tF_i|$ and $|\tF_j|$ on $Y$ have singular fibres of type
$I_0^*$ or $I_2^*$,
$G_i=2C_i+N_{i_1}+\dots +N_{i_4}$, respectively  $G_j=2C_j+N_{j_1}+\dots+
N_{j_4}$.\newline Then
$C_iC_j=0$ and the set $\{i_1,\dots, i_4\}\cap \{j_1,\dots, j_4\}$ consists of two
elements.
\end{lem}
\begin{proof} The curve $C_i$ is irreducible if  $G_i$ is of type $I_0^*$ and it is a
chain of 3 nodal curves if $G_i$ is of type $I_2^*$.  One has:
$4=\tF_i\tF_j=G_iG_j=2C_iG_j$, namely
$C_i(2C_j+N_{j_1}+\dots +N_{j_4})=2$.  We remark that $C_iN_{j_t}$ is equal to 1
 if $j_t\in \{i_1,\dots, i_4\}$, and it is equal to 0 otherwise. Since there are 7 of
the $N_i$, one has
$C_i(N_{j_1}+\dots +N_{j_4})>0$. So  either we have $C_iC_j=0$ and $\{i_1,\dots,
i_4\}\cap \{j_1,\dots, j_4\}$ consists of two elements, or $C_iC_j=-1$ and
$\{i_1,\dots, i_4\}=\{j_1,\dots, j_4\}$.

Assume by contradiction that we are in the second case. This implies in particular
that
$C_i$ and $C_j$ are not both irreducible. Assume that
$C_i$ is irreducible. Then   $C_j$ is a chain of 3 nodal curves
$C_j=A_1+N+A_2$ such that each of the ``end'' curves
$A_1$,
$A_2$  meets exactly two of the curves
$N_{i_1},\dots, N_{i_4}$ and the ``central'' curve $N$ is one of the $N_i$. In fact,
if $N$ were not one of the $N_i$, then the classes of $N_1,
\dots, N_7, N, A_1, A_2$, being independent,   would be a basis of
$H^2(Y,\Q)$, against the Index Theorem. Furthermore, from $C_iC_j=-1$ it follows that
$C_i$ is equal to
$A_1$ or $A_2$, a contradiction, since
$C_i$ meets all the curves $N_{i_1},\dots, N_{i_4}$. So we have $C_j=A_1+N+A_2$ as
above and, with an analogous notation,
$C_i=B_1+N'+B_2$, where
$N'$ is again  one   of the $N_i$.

Observe that  $\theta C_i\ge -1$ for every irreducible curve $\theta$. Since
$N'$ and $N$ are different from $N_{i_1},\dots, N_{i_4}$, the relations $N'\tF_j=0$,
$N\tF_i=0$ give:
$N'C_j=NC_i=0$ and
$-1=C_iC_j=(A_1+A_2)C_i$,  hence, say, $A_1C_i=-1$, $A_2C_i=0$. So we can assume that
$A_1=B_1$, while $A_2$ is disjoint from $C_i$ and $B_2$ is disjoint from  $C_j$. Say
that $A_1=B_1$ meets the curves $N_{i_1}$ and $N_{i_2}$. Then the connected divisor
$\Delta=N_{i_1}+N_{i_2}+A_1+N+N'$ is orthogonal to both $\tF_i$ and $\tF_j$, so its
support is contained in both in
$G_i$ and  $G_j$. If $N\ne N'$, then the intersection form on the components of
$\Delta$ is semidefinite, hence by Zariski's Lemma $\Delta$ is the support of both
$G_i$ and $G_j$, but this is impossible. Hence $N=N'$, but this contradicts the fact
that $A_2$ and $C_i$ are disjoint.
\end{proof}

\begin{lem}\label{listafibre}
\begin{enumerate}
\item The fibres with reducible support that  occur in the pencils $|\tF_i|$ can be
of  the following types:
$I_2$, $_2I_2$, $I_0^*$,  $I_2^*$;
\item each pencil $|\tF_i|$ has at least a fibre of type $I_0^*$ or $I_2^*$;
\item a fibre of type $I_0^*$ of $|\tF_i|$ contains 3 or 4 of the $N_i$, each with
multiplicity 1, and a fibre of type
$I_2^*$ contains 4 of the $N_i$ with multiplicity 1 and one with multiplicity 2.
\end{enumerate}
\end{lem}
\begin{proof} We recall first of all that the multiple fibres of an elliptic pencil
are of type $mI_k$, $k\ge 0$ (\cite{bpv}, Ch. V, \S 7) and that the multiple fibres
of  an elliptic pencil on  an Enriques surface are precisely two double fibres (\cite{bpv}, Ch. VIII). The  nodal curves
$N_1,\dots, N_7$ are contained in fibres of
$|\tF_i|$ for $i=1,2,3$. For every singular fibre
$F_s$ of
$|\tF_i|$, we denote by $r(F_s)$ the number of irreducible curves contained in
$F_s$ and different from $N_1,\dots, N_7$. Since the subspace orthogonal to the class
of
$\tF_i$ in $H^2(Y,\Q)$ has dimension 9, Zariski's Lemma implies that
$8+\sum_{F_s}(r(F_s)-1)\le 9$, namely $r(F_s)\le 2$ for every singular fibre
$F_s$ of $|\tF_i|$ and there is at most one singular fibre $F_s$ with
$r(F_s)=2$. This shows that the possible types are
$_mI_2$, $_mI_3$, $_mI_4$, $I_0^*$, $I_1^*$, $I_2^*$, $III$, $IV$ and that, except
possibly one,  the fibres with  reducible support are of type
$_mI_2$, $III$ or
$I_0^*$. On the other hand, we have $12=c_2(Y)=\sum_{F_s\ {\rm singular}}e(F_s)$,
hence the quantity $\sum_{F_s\ {\rm reducible}}e(F_s)$ is $\le 12$. Using this remark
and the fact that the  7 curves $N_1,\dots, N_7$ are contained in fibres of
$|\tF_i|$ it is easy to show that types $_mI_3$, $I_1^*$, $III$ and $IV$ cannot occur
and that the fibres cannot all be of type $_mI_2$ or $_mI_4$. This proves (ii).

Now assume that, say, $|\tF_1|$ has a fibre
$G_1$ of type $I_4$ (or
$_2I_4$). Then $G_1$ contains two of the $N_i$, hence we can write  the support of
$G_1$ as $C_1+C_2+N_1+N_2$, where $C_1C_2=0$, $C_iN_j=1$. Assume that $C_1\tE_2=0$.
Then the connected fundamental cycle $C_1+N_1+N_2$ is contained in a reducible fibre
$G_2$ of $|\tF_2|$. Since $C_1N_i=0$ for $i>2$, the fibre
$G_2$ is necessarily of type $I_4$ (or $_2I_4$). Since $r(G_1)=r(G_2)=2$ by the above discussion and
by (ii) it follows that both $|\tF_1|$ and $|\tF_2|$ have  fibres $G_1'$, respectively $G_2'$,
 of type
$I_0^*$  and that the nodal  curves appearing with multiplicity 1 in $G_1'$ and  $G_2'$ are a subset
of $\{N_3,\dots, N_7\}$,   contradicting Lemma
\ref{scemo}. This shows that the intersection numbers  $C_1\tE_2, C_1\tE_3,
C_2\tE_2,C_2\tE_3$  are all strictly  positive. Since
$4=(\tE_2+\tE_3)\tF_1=(\tE_2+\tE_3)G_1\ge C_1\tE_2+C_1\tE_3+ C_2\tE_2+C_2\tE_3$,
these numbers are all equal to 1. So the class of
$C_1-C_2$ is orthogonal to
$N_1,\dots, N_7, \tE_1,\tE_2,\tE_3$. Since the classes of $N_1,\dots, N_7,
\tE_1,\tE_2,\tE_3$ are a basis of
$H^2(Y,\Q)$, the class $C_1-C_2$ is numerically equivalent to 0. On the other hand,
we have
$(C_1-C_2)^2=-4$, a contradiction. This finishes the proof of (i).

 Statement (iii) follows by examining the admissible  types of fibres, recalling that
$r(F_s)\le 2$ for every singular fibre $F_s$.
\end{proof}

\begin{lem}\label{extranode} Assume that there exists a nodal curve $C\subset Y$ such that
$\tE_1C= N_1C=\dots =N_7C=0$. Then we have case 1) of Theorem \ref{B+fibres}.
\end{lem}
\begin{proof} By \cite[Lemma 4.2]{mp3}, the divisor $C+N_1+\dots+N_7$ is divisible by
2 in
$\Pic(Y)$. Hence
$C\tE_i$ is even for $i=1,2,3$. The curve $C$ is contained in a fibre of
$|\tF_1|$, hence
$C(\tE_2+\tE_3)=C\tB\le\tB\tF_1=4$.

 On the other hand, since $\tB+N_1+\dots +N_7$ is also even, the divisor $\tB+C$ is
even and so $(\tB+C)^2=4+2\tB C$ is divisible by 8. Hence we have $\tB C=2$. From
$2=\tB C=C\tE_2+C\tE_3$, it follows, say, $C\tE_2=0$, $C\tE_3=2$. Now, as in the
proof of Lemma \ref{listafibre}, we consider the contributions to $c_2(Y)$ and to the
Picard number of $Y$ of the various types of singular fibres. Since there are 8
disjoint nodal curves contained in the fibres of $|\tF_1|$ and
$|\tF_2|$, one sees  that the only possibility is that the fibres with singular
support of both  pencils are two fibres of type $I_0^*$ and that each fibre of type
$I_0^*$ contains four of the curves $N_1, \dots ,N_7, C$, each with multiplicity 1.
Recall that an elliptic pencil with 2 fibres of type $I_0^*$ on an Enriques surface
is isotrivial.

 By Lemma \ref{scemo} we can label the curves $N_i$ in such a way that the singular
fibres of $|\tF_1|$ are $N_1+N_2+N_3+N_4+2A_1$ and
$N_5+N_6+N_7+C+2A_2$ and the singular fibres of $|\tF_2|$ are $N_1+N_2+N_5+N_6+2B_1$
and
$N_3+N_4+N_7+C+2B_2$.

  Computing intersection numbers, one sees that
$\tE_3\sim_{num}
\tE_1+\tE_2-C$, namely
$\tE_3\equiv \tE_1+\tE_2-C$ or $\tE_3\equiv \tE_1+\tE_2-C+K_Y$. Since
$\tB+C$ is even, we conclude that $\tE_3\equiv
\tE_1+\tE_2-C$. Hence the system $|\tB|$ is equal to $|\tF_1+\tF_2-C|$ and it has degree
2 by \cite[Thm. 4.7.2]{codo}. Looking at the adjunction sequence  for $C$, one gets
$0=h^0(C+K_Y)=h^0(\tE_1+\tE_2-\tE_3')=0$, hence $|\tB+K_Y|$ is birational again by
\cite[Thm. 4.7.2]{codo}.

Now
$\tF_3\equiv
\tF_1+\tF_2-2C\equiv N_5+N_6+N_7+2A_2 +N_3+N_4+N_7+2B_2$, hence
$N_3+N_4+N_5+N_6+2(A_2+N_7+B_2)$ is a fibre of
$|\tF_3|$ of type
$I_2^*$. Thus Lemma
\ref{listafibre} and the formula $12=c_2(Y)=\sum_{F_s}e(F_s)$ imply that the
remaining fibres with singular support are two fibres of type
$I_2$ or
$_2I_2$.
\end{proof}

\begin{proof}[Proof of Theorem \ref{B+fibres}] By Lemma \ref{extranode}, we may
assume that every nodal curve disjoint from
$N_1, \dots, N_7$ satisfies $C\tE_i>0$ for $i=1,2,3$. Hence none of the pencils
$|\tF_i|$ can have two fibres of type
$I_0^*$.
 By Lemma
\ref{listafibre} and its proof one sees that in principle the possible configurations
of fibres with reducible  support are:
\begin{enumerate}
\item one fibre of type $I_0^*$ and three fibres of type  $I_2$ or
$_2I_2$;
\item one fibre of type $I_2^*$ and two fibres of type $I_2$ or
$_2I_2$.
\end{enumerate} Notice that in both  cases  all the  fibres with singular support are
reducible.

Assume that, say, $|\tF_1|$ has a fibre  $G_1$ of type $I_2^*$. By Lemma
\ref{listafibre}, (iii),
$G_1$ contains five of the $N_i$ and two more components  $A_2$ and $A_3$. Each of
the curves
$A_2$ and $A_3$ meets 3 of the $N_i$ and there is only one of the $N_i$ that
intersects both. We set $\lambda_2=2A_2+\sum_i(A_2N_i)N_i$ and
$\lambda_3=2A_3+\sum_i(A_3N_i)N_i$. One has
$\la_2^2=\la_3^2=-2$ and $G_1=\lambda_2+\lambda_3$. Since
$2=\tE_2G_1=2\tE_2(A_2+A_3)$, we may assume $\tE_2 A_2=1$, $\tE_2A_3=0$. Since by Corollary
\ref{Bample} the curves $N_1,
\dots, N_7$ are the only nodal curves orthogonal to $\tE_1, \tE_2, \tE_3$, we have
also $\tE_3A_3=1$, $\tE_3A_2=0$. The support of $\la_3$, being connected and
orthogonal to $\tE_2$,  is contained in a fibre $G_2$ of
$|\tF_2|$. Since $A_3$ meets precisely three of the $N_i$, by Lemma
\ref{listafibre} the fibre $G_2$ is also of type
$I_2^*$ and we can write as above $G_2=\la_1+\la_3$, where
$\la_1=2A_1+\sum_i(A_1N_i)N_i$, with $A_1$  a nodal curve different from the $N_i$
and such that
$A_1\tE_3=0$,
$A_1\tE_1=1$. Notice that the three nodal curves $A_1$, $A_2$ and $A_3$ are distinct.
The same argument shows that
$\la_1$ and
$\la_2$ are contained in fibres of
$|\tF_3|$ of type $I_2^*$. By the proof of Lemma \ref{listafibre} each pencil $|\tF_i|$ has at most
one fibre of type $I_2^*$, hence $\la_1$ and $\la_2$  are contained in the same fibre
$G_3$ and
$G_3=\la_1+\la_2$. Assume that the curve $N_i$ that appears with multiplicity $2$ in
$G_1$ and
$G_2$ is the same, say $N_7$. Then $N_7$ is a component of $\la_1$, $\la_2$ and
$\la_3$ and Lemma \ref{scemo} implies that, up to a  permutation of  $1,\dots, 6$ the
incidence relations between the curves
$N_1,\dots, N_7, A_1,A_2,A_3$ are given by the dual graph of  Figure
\ref{I22} below.

\begin{center}
\setlength{\unitlength}{0.5truecm}
\begin{picture}(16,15)\label{I22}
\put(5,2){\circle*{0.25}}
\put(11,2){\circle*{0.25}}
\put(5,4){\circle*{0.25}}
\put(11,4){\circle*{0.25}}
\put(3,4){\circle*{0.25}}
\put(13,4){\circle*{0.25}}
\put(8,7){\circle*{0.25}}
\put(8,11){\circle*{0.25}}
\put(7,13){\circle*{0.25}}
\put(9,13){\circle*{0.25}}
\thicklines
\put(5,2){\line(0,1){2}}
\put(11,2){\line(0,1){2}}
\put(8,7){\line(0,1){4}}
\put(3,4){\line(1,0){2}}
\put(11,4){\line(1,0){2}}
\put(5,4){\line(1,1){3}}
\put(11,4){\line(-1,1){3}}
\put(8,11){\line(-1,2){1}}
\put(8,11){\line(1,2){1}}
\put(7.5,5.5){$N_7$}
\put(5.3,3.5){$A_1$}
\put(9.7,3.5){$A_2$}
\put(8.5,10.5){$A_3$}
\put(4.5,0.8){$N_6$}
\put(10.5,0.8){$N_1$}
\put(1.5,3.8){$N_5$}
\put(13.5,3.8){$N_2$}
\put(5.5,12.8){$N_4$}
\put(9.3,12.8){$N_3$}
\end{picture}

Figure \ref{I22}
\end{center} The divisor $A_1+A_2+A_3+N_1+\dots+N_7$ is simply connected, hence its
inverse image in the K3 cover of $Y$ consists of two disjoint divisors isomorphic to
it. It is easy to check that the intersection matrix of the components of these
divisors is nondegenerate of type $(2, 18)$, but this contradicts the Index Theorem.
So, up to a permutation of the indices,
$G_1$ contains
$N_1$ with multiplicity 2, $G_2$ contains $N_2$ with multiplicity 2 and $G_3$
contains $N_3$ with multiplicity 2. Using Lemma \ref{scemo} again, one shows that, up
to a permutation of the indices, $N_7$ is not contained in
$G_1$,
$G_2$,
$G_3$ and the incidence relations between $A_1, A_2, A_3, N_1,\dots, N_6$ are given
by the dual graph in Figure \ref{I21}. Computing intersection numbers, one checks
that $\tB$ and the divisor
$\la_1+\la_2+\la_3$ are numerically equivalent. The argument used to  prove  Corollary \ref{corex2}
shows that they are not linearly equivalent, hence $\la_1+\la_2+\la_3\equiv \tB+K_Y$.
Now the system
$|\la_1+\la_2+\la_3|$ has degree 2 by Theorem 7.2 of \cite{cossec} and $|\tB|$ is
birational. This settles case 2).

We are left with the case in which each of the pencils $|\tF_i|$ has a fibre of type
$I_0^*$ and 3 fibres of type $I_2$ or $_2I_2$. We recall that  by \cite[Prop.
5.2.1]{cossec} and \cite[Thm. 4.7.2]{codo}   $|\tB|$ has degree 2   if
$|\tE_1+\tE_2-\tE_3|$ is nonempty and it is birational otherwise. Assume that there
is  $\Delta\in |\tE_1+\tE_2-\tE_3|$. Since
$\tE_1\De=\tE_2\De=0$ and $\tE_3\De=2$, all the components of $\De$ are nodal curves
contained in fibres of $|\tF_1|$ and $|\tF_2|$ and there is a component $\theta$ of
$\De$ with $\theta
\tE_3>0$, but this cannot happen because of the configuration of reducible fibres of
the pencils $|\tF_i|$. Since this argument is purely numerical it  shows also that
$|\tB+K_Y|$ is birational.
\end{proof}

\section{Enriques surfaces with 7 nodes: a general construction}\label{secconstr}

Here we describe a construction giving pairs $(\Si, B)$  as in the  set-up of \S
\ref{secex} and we prove that all such pairs  can be obtained that way.

\medskip

\begin{constr}\label{construction}{\em Consider the following automorphisms of the
projective line $\pp^1$:
$$(x_0, x_1)\stackrel{e_1}{\mapsto}(x_0,-x_1);\qquad
(x_0,x_1)\stackrel{e_2}{\mapsto}(x_1,x_0).$$

The subgroup $\Ga$ generated by $e_1$ and $e_2$ is isomorphic to $\Z_2^2$; we set
$e_3:=e_1+e_2$. The action of $\Ga$ can be lifted to  the line bundle
$\OO_{\pp^1}(2)$. Two such liftings differ by a character of $\Ga$,  hence for any chosen lifting a
basis of eigenvectors of
$H^0(\pp^1,\OO_{\pp^1}(2))$  is given by:
 $$s(x_0,x_1):=x_0^2+x_1^2,\quad d(x_0, x_1):=x_0^2-x_1^2,\quad p(x_0,x_1)=x_0x_1.$$ 

 Denote by $G$ the subgroup of automorphisms of
$\pp^1\times\pp^1\times\pp^1$ generated by the  elements:
 $$(e_1, e_1, 1),\ (e_1,1,e_1),\ (e_2,e_2,e_2).$$ The group $G$ is isomorphic to
$\Z_2^3$. We denote by
$G_0$ the subgroup of index 2 generated by $(e_1, e_1, 1)$ and
$(e_1,1,e_1)$.  The fixed locus of the nonzero elements of $G_0$ has dimension 1,
while the fixed locus of the elements of $G\setminus G_0$ has dimension 0.

Notice that, although the action of $\Ga$ on $\pp^1$ does not lift to a linear
representation on the space $H^0(\OO_{\pp^1}(1))$, the action of $G$ on
$\pp^1\times\pp^1\times\pp^1$ is induced by a linear representation on
$H^0(\OO_{\pp^1\times\pp^1\times\pp^1}(1,1,1))$. Hence it is possible to lift the
$G-$action to the line bundle $\OO_{\pp^1\times\pp^1\times\pp^1}(1,1,1)$ and,
compatibly, to all its multiples. Notice also that the possible $G-$actions on
$\OO_{\pp^1\times\pp^1\times\pp^1}(1,1,1)$ differ by a character of $G$, and thus 
they all induce the same action on $\OO_{\pp^1\times\pp^1\times\pp^1}(2,2,2)$.
Denoting the  homogeneous coordinates on
$\pp^1\times\pp^1\times\pp^1$ by  $x=(x_0,x_1)$, $y=(y_0,y_1)$,
$z=(z_0,z_1)$, under this action  the space
$H^0(\OO_{\pp^1\times\pp^1\times\pp^1}(2,2,2))$ decomposes into eigenspaces as
follows:
$$T_0:=<s(x)s(y)s(z), s(x)d(y)d(z), d(x)s(y)d(z),$$ $$ d(x)d(y)s(z), p(x)p(y)p(z)>;$$
$$T_1:=<s(x)s(y)d(z), s(x)d(y)s(z), d(x)s(y)s(z), d(x)d(y)d(z)>;$$
$$T_2:=<s(x)s(y)p(z), p(x)p(y)s(z), d(x)d(y)p(z)>;$$
$$T_3:=<s(x)p(y)p(z), p(x)s(y)s(z), p(x)d(y)d(z)>;$$
$$T_4:=<p(x)s(y)d(z), p(x)d(y)s(z), d(x)p(y)p(z)>;$$
$$T_5:=<s(x)p(y)s(z), p(x)s(y)p(z), d(x)p(y)d(z)>;$$
$$T_6:=<s(x)d(y)p(z), p(x)p(y)d(z), d(x)s(y)p(z)>;$$
$$T_7:=<s(x)p(y)d(z), p(x)d(y)p(z), d(x)p(y)s(z)>.$$

The subspace $T_0$ corresponds to the trivial character and the subspace
$T_1$ corresponds to the  character orthogonal to $G_0$. The system $|T_0|$ is base
point free, while for
$i>0$ the base locus of the  system
$|T_i|$ is nonempty and it contains the fixed locus of some   element of
$G\setminus G_0$.

Let $Z\in |T_0|$ be general. By Bertini's theorem $Z$ is a smooth surface. The
elements $(e_1, 1,e_1)$, $(1,e_1, e_1)$ and $(e_1, e_1, 1)$ act on $Z$ fixing 8
points each and the remaining nonzero elements of $G$ act freely on $Z$. By the
adjunction formula, $Z$ is a K3 surface, hence the quotient surface $Z/G$ is a nodal
Enriques surface with 6 nodes. For $i=1,2,3$, the   projection on the $i-$th factor
$\pp^1\times\pp^1\times\pp^1\to\pp^1$ induces an elliptic pencil on $Z$, which in turn
gives an elliptic pencil $|F_i|=|2E_i|$ on $Z/G$. A standard argument shows that for
$i\ne j$ one has $E_iE_j=1$.

Assume now that $Z$ has 8 nodes that form a $G-$orbit and no other singularities. Then
the quotient surface $Z/G$ has an extra node, which is the image of the 8 nodes of
$G$. By Lemma \ref{even} either  $E_1+E_2+E_3$ or $E_1+E_2+E_3+K_{Z/G}$ is a divisor
$B$ as in the set-up of \S \ref{secex}.

Set $\Si:=Z/G$ and denote by  $\pi\colon Z\to \Si$ the quotient map. Notice that
$\pi^*B\iso \pi^*(B+K_{\Si})$ is isomorphic to
$\OO_Z(2,2,2)$ and
$\pi^*H^0(\Si,B)$ and $\pi^*H^0(\Si, B+K_{\Si})$ are eigenspaces of
$H^0(Z,
\OO_Z(2,2,2))$. So,  considering the dimensions,  they correspond to the restrictions
to
$Z$ of $T_0$ and
$T_1$. We will show later (Lemma \ref{wiw}) that the restriction of $T_0$ is equal to
$\pi^*H^0(\Si, B+K_{\Si})$ and the restriction of $T_1$ is equal to
$\pi^*H^0(\Si,B)$.}
\end{constr}

The fact that Construction \ref{construction}   can actually be performed, namely
that there exists $Z$ as required, is a consequence  of  the following Theorem and of
the examples given in \S \ref{secex}.
\begin{thm}\label{structure}
 Let $(\Si,B)$ be a pair as in the set-up of \S \ref{secex}. Then $(\Si,B)$ can be
obtained from Construction \ref{construction}.
\end{thm}
\begin{proof} Let $V$ be the code associated with the nodes of $\Si$, which is
isomorphic to $\Z_2^2$
 by Proposition \ref{dimV}. By \cite{nodes}, Prop. 2.1 and Remark 2, there is a
Galois cover
$\pi_0\colon Z_0\to\Si$ with Galois group $\Hom(V,\C^*)\iso \Z_2^2$ branched on the 6
nodes of $\Si$ that appear in $V$.
The map $Z_0\to\Si$ can be factorized as $Z\to Z_1\to \Si$, where both maps are
 double covers branched on a set of  4 nodes. By \cite[Prop.
3.1]{mp3},
$Z_1$ is a nodal Enriques surface with 6 nodes,  hence, by ib.,  $Z_0$ is  an
Enriques surface with 4 nodes.  Let $K\to \Si$ be the
K3-cover of $\Si$ and consider the following cartesian diagram:
 \begin{equation}\label{dia1}
\begin{CD} Z@>>>Z_0\\ @Vp VV@V\pi_0 VV\\ K@>>>\Sigma
\end{CD}
\end{equation} The surface $Z$ is a K3 surface with 8 nodes and the map $Z\to Z_0$ is
the K3-cover of
$Z_0$. The composite map $\pi\colon Z\to\Si$ is a Galois cover with Galois group
isomorphic to $\Z_2^3$.
 Notice (cf. \cite{nodes}, proof of Prop. 2.1 and Remark 2) that, although the cover
$\pi_0\colon Z_0\to\Si$ is not uniquely determined (in fact there are four different
possibilities), the cover
$\pi\colon Z\to \Si$ does not depend on the choice of $Z_0$.
\smallskip

For the reader's convenience the proof is broken into steps.

\noindent{\bf Step 1:} {\em For $i=1,2,3$ there exist elliptic  pencils
$|C_i|$ on $Z$ such that
$\pi^*F_i\equiv 4C_i$.}

Since $\pi$ is unramified in codimension 1, if $F_i\in |F_i|$ is  general then
$\pi^*F_i$ is a disjoint union of linearly equivalent elliptic curves. Hence, to
prove the statement it is enough to show that $\pi^*F_i$ has 4 connected components.
Let $\tpi\colon \tZ\to Y$ the Galois cover obtained from $\pi$ by taking base change
with the minimal desingularization
$Y\to\Si$ and, as usual, denote by $|\tF_i|$ the elliptic  pencil of $Y$ induced by
$|F_i|$.
 By Theorem \ref{B+fibres} the pencil $|\tF_i|$ has a fibre  $G_i$ of type
$I_0^*$ or
$I_2^*$. We write $G_i=2A_i+N^i_1+\dots+  N^i_4$. By the results of \S
\ref{seccode}, the nonzero elements of $V$ correspond to the  even sets $N^i_1+\dots
+ N^i_4$,
$i=1,2,3$. So by the definition of
$Z$ and $\tZ$ (cf. also \cite[\S 2]{nodes}) we have the following formula:
$$\tpi_*\OO_{\tZ}=\OO_Y\oplus\-K_Y\oplus(\oplus_{i=1,2,3}
\OO_Y(A_i-\tE_i))\oplus (\oplus_{i=1,2,3}
\OO_Y(A_i-\tE_i')).$$

The restriction  of the line bundles $\OO_Y(-K_Y)$,
$\OO_Y(A_i-\tE_i)$ and
$\OO_Y(A_i-\tE_i')$ to a general $\tF_i$ is trivial, hence $\tpi^*\tF_i$ has at least
4 connected components. So for $i=1,2,3$ we can write $\tpi^*F_i\equiv m_i\tC_i$ 
where $\tC_i$ is a smooth connected elliptic curve and $m_i=4$ or $m_i=8$. Notice
that $\tC_i\tC_j\ge 2$ for
$i\ne j$, since otherwise the product of the pencils $|\tC_i|$ and $|\tC_j|$ would
give  a birational map $\tZ\to \pp^1\times
\pp^1$. On the other hand, for $i\ne j$ we have
$32=\tpi^*\tF_i\tpi^*\tF_j=m_im_j\tC_i\tC_j\ge 16
\tC_i\tC_j$, hence
$m_i=m_j=4$ and $\tC_i\tC_j=2$. Finally,  the pencils $|\tC_i|$ induce pencils
$|C_i|$ on $Z$ such that $\pi^*F_i\equiv 4C_i$.
\medskip

\noindent{\bf Step 2:} {\em The product of the pencils $|C_1|$, $|C_2|$ and
$|C_3|$ defines an embedding  $\psi\colon Z\into \pp^1\times \pp^1\times \pp^1$ such
that
$\psi(Z)$  is a divisor of type
$(2,2,2)$.}
 We remark that $D=C_1+C_2+C_3$ is ample on $Z$, since by Step 1  $2D=\pi^*B$ and $B$
is ample on $\Si$ by Corollay \ref{Bample}.

We let $\epsi\colon W\to Z$ be the minimal desingularization, we denote by
$N_1',\dots ,N_8'\/$ the exceptional curves of $\epsi$ and we set
$\tD:=\epsi^*D$, $\tC_i:=\epsi^*C_i$. The divisor
$\tD$ is nef and big and the
$N'_i$ are the only irreducible curves that have intersection equal to 0 with $\tD$.
Since
$\tD^2=12$, by Reider's Theorem and by the fact that the intersection form on a K3
surface is even, if  two points $x,y\in W$ are not separated by $|\tD|$, then there
exists an effective connected divisor
$A\ni x,y$ such that either: a) $A^2=-2$, $A\tD =0$, or b) $A^2=0$, $A\tD=1$ or c)
$A^2=0$,
$A\tD=2$.
 Possibility a)  corresponds to the case when both $x$ and $y$ belong to one of the
curves
$N'_i$. Recall that the pencils $|\tC_i|$ have no multiple fibres, because the double
fibres of the pencils $|\tF_i|$ disappear when one takes the K3 cover (actually, it
is not hard to prove that any elliptic fibration on a K3 has no multiple fibres). In 
case b) one would have, say,
$A\tC_1=A\tC_2=0$, namely
$A$ would be a fibre of both
$|\tC_1|$ and $|\tC_2|$, which is impossible. If $A^2=0$ and $A\tD =2$, then we have,
say, $A\tC_1=0$ and $A$ is a fibre of $|\tC_1|$. But in this case $A\tD=4$, a
contradiction.

The above discussion shows that the map $\psi$ is one-to-one onto its image and that
the differential of
$\psi$ at every smooth point of $\Si$ is nonsingular. In particular the image of
$\psi$ is   an hypersurface with at most isolated singularities, hence it is  
normal. It follows that
$\psi$ is an isomorphism. The fact that the image is a divisor of type $(2,2,2)$ is a
consequence of the fact that
$C_iC_j=2$ if $i\ne j$.
\smallskip

\noindent{\bf Step 3:} {\em There exist coordinates on
$\pp^1\times\pp^1\times\pp^1$ such that the surface
 $\psi(Z)$ is an element of the linear system $|T_0|$ defined in Construction
\ref{construction} and the action of the Galois group of $Z\to\Si$ coincides with
the  group action  defined there.}

Denote by $G$ the Galois group of $\pi$. By the definition of the map $\psi$, the
three copies of $\pp^1$ in
$\pp^1\times\pp^1\times\pp^1$ can be naturally identified with (the dual of)
 the  linear systems $|C_1|$,
$|C_2|$ and
$|C_3|$,  hence $G$ acts on $\pp^1\times\pp^1\times\pp^1$ and the embedding
$\psi\colon Z\to
\pp^1\times\pp^1\times\pp^1$ is $G-$equivariant with respect to the given actions. We
have seen in Step 1 that for every
$i=1,2,3$ there is a nonzero $g_i\in G$ such that $g_i$ acts trivially on
$|C_i|$. Since the fixed locus of $g_i$ on $|C_1|\times |C_2|\times |C_3|$ has positive dimension and
$Z$ is ample,  
$g_i$ has fixed points on $Z$.  Since by construction the cover
$\pi\colon Z\to\Si$ factorizes through the K3 cover $K\to\Si$, it follows   that
$g_1$, $g_2$, $g_3$ do not generate $G$. On the other hand
$g_i$ must act non trivially on $|C_j|$ for $j\ne i$, since otherwise the fixed locus
of $g_i$ on
$\pp^1\times\pp^1\times\pp^1$ would be a divisor and $g_i$ would fix a curve of $Z$
pointwise. Hence $G_0:=\{1, g_1, g_2,g_3\}$ is a subgroup of $G$ isomorphic to
$\Z_2^2$. Fix $h\in G\setminus G_0$. For every $i$ we can choose homogeneous
coordinates on $\pp^1=|C_i|$ such that, using the notation of Construction
\ref{construction}, the nonzero element of $G_0/g_i$ acts as $e_1$ and $h$ acts as
$e_2$. With respect to these coordinates we have: $g_1=(1,e_1, e_1)$, $g_2=(e_1, 1,
e_1)$,
$g_3=(e_1, e_1,1)$, $h=(e_2, e_2, e_2)$, namely the $G-$action on
$\pp^1\times\pp^1\times\pp^1$ is the same as in Construction \ref{construction} and
the surface  $Z$, being $G-$invariant, belongs to one of the linear systems $|T_i|$,
$i=0,\dots, 7$. In addition, each of the nonzero elements of
$G_0$ fixes 8 points of $Z$ and the elements of $G\setminus G_0$ act freely on $Z$.
This is the same as saying that $Z$ is in general position with respect to the fixed
loci of all the elements, hence, as we have remarked in Construction
\ref{construction}, $Z$ must be  an element of $|T_0|$.

\end{proof}

\section{Enriques surfaces with 7 nodes: a parametrization}\label{secmod}

The aim of this section is to construct a quasi-projective  variety parametrizing the
isomorphism classes of pairs
$(\Si,B)$ as in the set up of \S \ref{secex}  and to study  the geometry of this
space. In addition,  we show the existence of a tautological family on a finite Galois
cover of the parametre space. This tautological family admits a section and  a
simultaneous resolution.

 These results are used  in the next section to describe the subset of the moduli
space of surfaces  with $p_g=0$ and $K^2=3$ consisting of the surfaces $S$ that have
an involution
$\si$  such that: 1) the quotient surface $S/\si$ is birational to an Enriques
surface; 2) the bicanonical map
$\fie$ of $S$ is composed with $\si$.
\smallskip

 We use all the notation from  the previous sections. For a pair $(\Si, B)$ as in the
set--up  of \S \ref{secex}  we denote   as usual by $\eta\colon Y\to
\Si$ the minimal desingularization and by $N_1,\dots , N_7$ the exceptional curves of
$\eta$. In addition we assume that $N_7$ is the nodal curve that does not appear in
the code $V$ associated to
$N_1,\dots, N_7$ (cf. \S \ref{seccode}).

Denote by $\cN$ the subset of $|T_0|$ consisting of the surfaces $Z$ that satisfy the
following conditions:

a) $Z$ is in general position with respect to the fixed loci of the elements of $G$;

b) $Z$ has 8 nodes that form a $G-$orbit and no other singularities. The set $\cN$ is
clearly open in the set of singular surfaces of $|T_0|$ and it is nonempty by the
results of the previous section, hence it is a quasi-projective variety of dimension
3. We denote by $\cI\subset (\pp^1\times\pp^1\times \pp^1)\times \cN$ the incidence
variety, consisting of the pairs $(P,Z)$ such that $P$ is a singular point of $Z$ and
we denote by $p_1$,
$p_2$ the projections of $\cI$ onto the two factors.   There is a natural $G-$action
on $\cI$, which is free by the definition of $\cN$, and the map
$p_2\colon\cI\to\cN$ is the quotient map with respect to this $G-$action. 

The first goal of this section is to study the geometry of $\cN$. We have the
following: 
\begin{thm}\label{Ngeom} The  variety  $\cN$ is smooth, irreducible of dimension 3
and  unirational. 
\end{thm}

The proof that $\cN$ is smooth is completely elementary (cf. Lemma \ref{Nsmooth}
below),  but proving the irreducibility requires a  series of intermediate results.

\begin{lem}\label{Nsmooth} The variety $\cN$ is smooth of dimension 3.
\end{lem}
\begin{proof} Since the incidence variety $\cI$ is a topological covering of $\cN$,
it is enough to prove that
$\cI$ is smooth of dimension 3. This can be easily seen by means of  a local
computation, using the fact that the  linear system $|T_0|$ has no base points and
the fact that for a pair 
$(P,Z)$  in
$\cI$ the point
$P$ is an ordinary double point of $Z$.
\end{proof}
\begin{lem} \label{wiw} Assume that the pair  $(\Si, B)$ is obtained from
$Z\in \cN$ using Construction
\ref{construction} and let $\pi\colon Z\to\Si$ be the quotient map. Then:
$$\pi^*H^0(\Si, B)=T_1|_Z;\qquad \pi^*H^0(\Si, B+K_{\Si})=T_0|_Z.$$
\end{lem}
\begin{proof} We  have already remarked in Construction \ref{construction} that 
$H^0(\Si, B)$ and $H^0(\Si, B+K_{\Si})$ pull back to $T_0|_Z$ and $T_1|_Z$. So we
only need to decide which is which.

Let $s\in T_0|_Z$ be general,  let $D$ be the divisor of zeros of $s$ and let
$\overline{D}$ be the image of
$D$ in
$\Si$. The divisor
$\overline{D}$ is smooth and it is numerically equivalent to  $B$. Let $f\colon X\to
Z$ be the double cover branched on $D$. Denote by $L$ the total space of the line
bundle $\OO_Z(1,1,1)$, by $p\colon L\to Z$ the projection and by $z$ the tautological
section of $p^*\OO_Z(1,1,1)$. Then $X$ is isomorphic to the hypersurface
$\{z^2-p^*s=0\}\subset L$ and the $G-$action on $L$ (cf. Construction
\ref{construction}) preserves
$X$. Hence the
$G-$action on $Z$ lifts to $X$ and we have a commutative diagram:
\begin{equation}\label{dia1}
\begin{CD} X@>>>Z\\ @Vq VV@V\pi VV\\
\overline{X} @>>>\Sigma
\end{CD}
\end{equation} where $q\colon X\to\overline{X}:=X/G$ is the quotient map. By
commutativity of the diagram, the map $\overline{X}\to\Si$ is a double cover branched
on $\overline{D}$ and on a subset of the 6 nodes of
$\Si$ that are the images of  the fixed points of the $G-$action on $Z$. As
before, let $\eta\colon Y\to\Si\/$\ be  the minimal desingularization of $\Si$.
Set
$\tD=\eta^*\overline{D}$ and denote by
$N_1,\dots, N_k$ the nodal curves of $Y$ corresponding to the nodes of $\Si$ where
$\overline{X}\to
\Si$ ramifies.   The class of
$\tD+N_1+\dots+N_k$ is divisible by 2 in $\Pic(X)$, hence its self-intersection,
which is  equal to $6-2k$, is  divisible by 8. Since $k\le 6$, it follows that $k=3$.
If  $\tD$ were linearly equivalent to $\tB$, then
$N_4+N_5+N_6+N_7$ would be  divisible by 2 in $\Pic(Y)$, contradicting the fact that
$N_7$ does not appear in the code $V$ associated with the curves $N_1,\dots,  N_7$.
So we must have $\tD\equiv
\tB+K_Y$, and thus
$T_0|_Z=\pi^*H^0(\Si, B+K_{\Si})$.
\end{proof}
\begin{lem}
 The  curve $N_7$ is not contained in a double fibre  of $|\tF_i|$, for
$i=1,2,3$.
\end{lem}
\begin{proof} Assume by contradiction that $N_7$ is contained in a double fibre of,
say,
$|\tF_1|$.
 Then, by Theorem \ref{B+fibres},  $N_7$ is contained in a fibre $2A$ of
$|\tF_1|$ with $A$ of type
$I_2$. The cover $\tpi\colon \tZ\to Y$ obtained from  $\pi\colon Z\to\Si$ by taking
base change with $\eta\colon Y\to\Si$ is \'etale over
$A$. More precisely, by Step 1 of the proof of Theorem \ref{structure} the divisor
$\tpi\inv(A)$ is the disjoint union of 2  connected curves, each mapping to $A$ with
Galois group $\Z_2^2$, but this is impossible since the fundamental group of $A$ is
cyclic.
\end{proof}
\begin{lem}\label{bir} Assume that $|K_{\Si}+B|$ is birational and let
$\psi\colon
\Si\to\pp^3$ be the corresponding morphism. Set $\Si':=\Si\setminus( E_1\cup
E_1'\cup\dots \cup E_3')$.

Then the restricted map  $\psi|_{\Si'}\colon \Si'\to\psi(\Si')$ is an isomorphism.
\end{lem}
\begin{proof} The map $\psi$ is a morphism onto a sextic of $\pp^3$. The divisor 
 $B$ is ample by Corollary \ref{Bample}, hence $K_{\Si}+B$ is also ample and
$\psi\colon \Si\to\psi(\Si)$ is the normalization map.  For
$i=1,2,3$,  the supports $E_i$, $E_i'$  of the double fibres
 of
$|F_i|$ are mapped $2-$to$-1$ onto distinct lines $L_i$,  $L_i'$ which are double for
$\psi(\Si)$. The general curve of $|K_{\Si}+B|$ is smooth of genus 4, hence the
general section
$C$  of
$\psi(\Si)$ has geometric  genus 4. Since $C$ has arithmetic genus 10 and it has at
least 6 singular points $C\cap L_1, \dots ,C\cap L_3'$, it follows that
$L_1,\dots,L_3'$ are the only $1-$dimensional components of the singular locus  of 
$\psi(\Si)$. Since
$K_{\Si}+B\equiv E_1+E_2+E_3'\equiv\dots
\equiv E_1'+E_2'+E_3'$,   the inverse image of $\psi(\Si)\setminus (L_1\cup\dots \cup
L_3')$ is $\Si'$.  The surface $\psi(\Si')=\psi(\Si)\setminus (L_1\cup\dots
\cup L_3')$ is normal, since it is an hypersurface and it is smooth in codimension 1.
It follows that  the map $\psi|_{\Si'}\colon
\Si'\to\psi(\Si')$ is an isomorphism.
\end{proof} We denote by $\cN_0\subset \cN$ the set of surfaces $Z$ such that $|T_0|$
induces a birational map $Z/G\to\pp^3$. By Lemma \ref{wiw} $Z\in \cN_0$ if and only
if the system $|K_{\Si}+B|$ is birational, where  $(\Si,B)$ is the pair obtained
from
$Z$ by  Construction \ref{construction}. The set
$\cN_0$ is open in $\cN$.

\begin{prop}\label{dense} The set $\cN_0$ is dense in $\cN$.
\end{prop}
\begin{proof} Since the proof is lengthy,  we describe first  the underlying idea,
which  is instead  quite simple. 

Let $Z\in\cN$ be a point.  Denote by 
$(\Si, B)$ the pair obtained from $Z$ by Construction \ref{construction} and denote by
$Z\to
\Si$ the
$G-$cover defined in the proof of Theorem \ref{structure}. Let $Y$ be the minimal
desingularization of $\Si$, let $N_1,\dots ,N_7$ be the corresponding nodal curves on
$Y$   and let
$\tZ\to Y$ be the $G-$cover obtained from $\pi\colon Z\to\Si$ by taking base change with
$\eta$.    As in Example \ref{ex3},  consider the restriction
$\cY_1\to (U_1,0)$  of the Kuranishi family of $Y$ to the subset where the classes of
$N_1,\dots, N_7$ stay effective. We show by standard arguments that one can construct
a $G-$cover $\widetilde{\cZ}\to \cY_1$ such that the fibre over 
$0\in U_1$ of the induced family of surfaces  
$\widetilde{\cZ}\to U_1$ is $\tZ$. Then, using the theory of \cite{bw} and a criterion for the
birationality of a linear system of degree 6 on an Enriques surface, we show that for
general $t\in U_1$ the surface $Z_t$ obtained by contracting the inverse images of the $N_i$  in
the  fibre
$\tZ_t$ of
$\widetilde{\cZ}$ is an element of 
$\cN_0$.
\smallskip

We start by constructing the family $\widetilde{\cZ}\to U_1$. 
 We recall that the base of the  Kuranishi family of $Y$ is smooth of dimension 10 by
\cite[Thm. VIII.19.3]{bpv}. Thus   the set 
$U_1$ is smooth of dimension 3 by
\cite[Thm. 2.14]{bw}.  Arguing as in Example \ref{ex3}, one sees that, after possibly
shrinking
$U_1$,  the building data of the $G-$cover
$\tZ\to Y$  can be extended to the total space $\cY_1$. Hence we have a
$G-$cover $\widetilde{\Pi}\colon \widetilde{\cZ}\to\cY_1$ that specializes to
$\tZ\to Y$ at the point
$0\in U_1$. The induced map $q\colon \widetilde{\cZ}\to U_1$ has smooth fibres. 

For
$i=1,2,3$, we  denote   by
$|\tC_i|$ the pull back to
$\tZ$ of the elliptic  pencil $|C_i|$ defined in Step 1 of the proof of Theorem
\ref{structure}. We claim that $\tC_i$ can be extended to a line bundle $\cC_i$ on
$\widetilde{\cZ}$ for
$i=1,2,3$. After possibly shrinking $U_1$  again, we may assume that $U_1$ is
contractible and that $\widetilde{\cZ}\to U_1$ is diffeomorphic to the product
family  $\tZ\times U_1$. Hence the inclusion $\tZ\to \widetilde{\cZ}$ induces an
 isomorphism $H^2(\widetilde{\cZ},\Z)\stackrel{\sim}{\to}H^2(\tZ,\Z)$. Recall that a
cohomology  class comes from a (holomorphic) line bundle on $\widetilde{\cZ}$ if and
only if it goes to zero under the map
$ob\colon H^2(\widetilde{\cZ},\Z)\to H^2(\widetilde{\cZ}, \OO_{\widetilde{\cZ}})$
induced by the exponential sequence. Arguing as in Example \ref{ex3}, one shows that there
exists a line bundle $\widetilde{\cF}_i$ on $\cY_1$ that restricts to $\tF_i$ on $Y$,
hence
$\widetilde{\Pi}^*\widetilde{\cF}_i$ induces the class of $\pi^*(\tF_i)$ and
$ob(\pi^*\tF_i)=0$. Since we have
$4\tC_i\equiv
\pi^*\tF_i$ by definition,  $ob(\tC_i)=\frac{1}{4}ob(\pi^*\tF_i)=0$ and the claim is
proven.

 Notice that the  line bundle representing  a given cohomology class is unique up to
isomorphism,  since
$H^1(\widetilde{\cZ},\OO_{\widetilde{\cZ}})=0$, as one  can check  by using the Leray
spectral sequence. The cohomology class of $\cC_i$ is $G-$invariant, since the class
of
$4\tC_i=\pi^*\tF_i$ is invariant and $H^2(\tZ,\Z)$ has no torsion. Thus  for every
$g\in G$ there is an isomorphism $g^*\cC_i\cong \cC_i$. Using this isomorphism,
every  $g\in G$ can be lifted to an automorphism of $\cC_i$.
 Hence we have a short exact sequence:
$$0\to\C^*\to  G_i\to G\to 0,$$ where $G_i$ is the group of automorphisms of $\cC_i$
that lift an element of
$G$. Since $h^0(\tZ,C_i)=2$ and $h^j(\tZ,C_i)=0$ for $j>0$, we may assume by the
semicontinuity theorem that
 $h^0(\tZ_t,\cC_i|_{\tZ_t})=2$ for all $t\in U_1$.  Hence for $i=1,2,3$,  the sheaf
$q_*\cC_i$ is a rank 2 vector bundle on $U_1$. We set $\cP_i:=\Proj(q_*\cC_i)$. By
the above exact sequence, $G$ acts on $\cP_i$ and  the natural map
$\Psi\colon\widetilde{\cZ}\to \cP_1\times \cP_2\times \cP_3$ is $G-$equivariant. We
remark that the restriction of $\Psi$ to $t=0$ is the composition of the map
$\tZ\to Z$ with the map $\psi$ defined in Step 2  of the proof of Theorem
\ref{structure}.
 The image of
$\widetilde{\cZ}$ is a divisor relatively of type $(2,2,2)$. Fix homogeneous
coordinates on  $|C_1|\times |C_2|\times |C_3|$ as in Step 3 of the
proof of Theorem \ref{structure}, namely such that the $G-$action is the one
described in  Construction \ref{construction}. By Step 3 of the proof of Theorem
\ref{structure}, in these coordinates 
$\psi(Z)$ is an element of
$|T_0|$.

We  want to construct a local trivialization
$\cP_1\times \cP_2\times\cP_3\stackrel{\sim}{\to}U_1\times (\pp^1\times
\pp^1\times \pp^1)$   that coincides for $t=0$ with the chosen isomorphism
$|C_1|\times |C_2|\times |C_3|\simeq \pp^1\times\pp^1\times\pp^1$ and
such that the
$G-$action on
$\cP_1\times \cP_2\times \cP_3$ corresponds fibrewise to the $G-$action on 
$\pp^1\times\pp^1\times \pp^1$  defined in Construction
\ref{construction}. We do this by finding trivializations of
$\cP_i$,
$i=1,2,3$, such that $\Ga=<e_1,e_2>$ acts as in Construction
\ref{construction}.    Start with any local trivialization and then  choose
homogeneous coordinates on
$\pp^1$  such that  the fixed locus of $e_1$ is $x_0x_1=0$ (this is possible since we
are working  in the analytic category). In these  coordinates the action of
$e_1$ is represented by  $(x_0, x_1)\mapsto (x_0,-x_1)$. Now $e_1$ exchanges the fixed
points of
$e_2$, which therefore have homogeneous  coordinates $(a(t), b(t))$, $(a(t), -b(t))$,
with $a$, $b$ nowhere vanishing functions on $U_1$ such that $a(0)=b(0)$.  Changing
coordinates to
$x'_0=b(t)x_0$,
$x_1'=a(t)x_1$, the fixed locus of $e_2$ is defined by $(x_0')^2-(x_1')^2=0$ and the
action of $e_2$ is induced  by $(x_0', x_1')\mapsto (x_1', x_0')$. Finally,  using 
an automorphism of $\pp^1$ that induces the identity on the group $\Ga$, one can
obtain that for $t=0$ the coordinates agree with the chosen ones. Now
$\Psi(\widetilde{\cZ})$ is defined inside $U_1\times (\pp^1\times\pp^1\times
\pp^1)$ by an equation $F(t, x,y,z)$ such that for fixed  $\bar{t}\in U_1$
$F(\bar{t},x,y,z)$ is a $G-$invariant  homogeneous polynomial of type $(2,2,2)$ and
$F(0,x,y, z)=0$ is the equation of $\psi(Z)$. Furthermore, the proof of Theorem
\ref{structure} implies that for every
$t\in U_1$ the surface
$F(t, x,y,z)=0$ is an element of $\cN$.

Assume now that $Z\notin \cN_0$.  By Lemma \ref{wiw}, the system $|B+K_{\Si}|$ is not
birational. The configuration of reducible fibres of the pencils $|\tF_i|$ on $Y$ is
described in Theorem \ref{B+fibres}, 2). By the proof of Theorem \ref{B+fibres},  the
components of the reducible fibres of the elliptic pencils can be labelled in such a
way that their dual diagram is the one given in Figure
\ref{I21} and $F_iA_j=2\delta_{i,j}$. In the notation of the proof of  Theorem
\ref{B+fibres}, we have $\la_1+\la_2+\la_3\equiv \tB+K_Y$ and $\la_1+\la_2\equiv
\tF_3$. Then we have
$N_1+N_2+N_6+ 2A_3=\la_3\equiv \tB+K_Y-\tF_3\equiv 
\tE_1+\tE_2-
\tE_3+K_Y$. We denote by $\widetilde{\cE}_i$  the line bundle on $\cY_1$ that extends
$\tE_i$, $i=1,2,3$.  By
\cite[Prop. 5.2.1]{cossec} and 
\cite[Thm. 4.7.2]{codo}, the system $|K_{Y_t}+\tB_t|$ is not birational if and only if
$\widetilde{\cE}_{1,t}+\widetilde{\cE}_{2,t}-\widetilde{\cE}_{3,t}+K_{Y_t}$ is
effective. Again  by
\cite[Thm. 3.7]{bw}, the subset
$U_2$ of the base $U$ of the Kuranishi family of $Y$ where
$\widetilde{\cE}_{1,t}+\widetilde{\cE}_{2,t}-\widetilde{\cE}_{3,t}+K_{Y_t}$,
$\cN_{3,t}$,
$\cN_{4,t}$, $\cN_{5,t}$,
$\cN_{7,t}$ are effective is smooth of dimension 2 in a neighbourhood of  $t=0$. Hence
$U_1\setminus U_2$ is nonempty and for $t\in U_1\setminus U_2$ the surface
$F(t,x,y,z)=0$ is an element of
$\cN_0$. This proves that $\cN_0$ is dense in $\cN$.
\end{proof}
\begin{proof}[Proof of Theorem \ref{Ngeom}] The variety $\cN$ is smooth of dimension
3 by Lemma \ref{Nsmooth}. By Lemma \ref{dense}, to complete the proof it suffices to
show  that $\cN_0$ is irreducible
 and unirational.

Denote by $\cI_0$ the restriction to $\cN_0$ of the incidence variety $\cI$ and denote
again by 
$p_1\colon
\cI_0\to\pp^1\times \pp^1\times\pp^1$ and $p_2\colon \cI_0\to \cN_0$ the projections.
The map
$p_2$ is an
\'etale $G-$cover by construction. We prove the theorem  by showing that
$p_1$  is injective.

Let $P\in p_1(\cI_0)$ be a point.  The fibre $p_1\inv(P)$ is a (nonempty) open subset
of the linear subsystem  of $|T_0|$ consisting of the surfaces which are singular at
$P$. Assume that $p_1\inv(P)$ has positive dimension and let $Z$, $Z'$ be two
distinct surfaces in $p_1\inv(P)$, let $(\Si,B)$ be the pair obtained from $Z$ by 
Construction \ref{construction} and let $\eta\colon Y\to \Si$ be the minimal
resolution. Recall that the image of $P$ is the node of $\Si$ corresponding to
the nodal curve $N_7$ of $Y$ that does not appear in the code $V$. Denote by
$D$ the divisor on
$\Si=Z/G$ induced by the restriction of
$Z'$ to
$Z$. By Lemma \ref{wiw}, $D\equiv B+K_{\Si}$.   It is not difficult to check that the
pull back
$\tD$ of
$D$ to
$Y$ vanishes on $N_7$ of order at least 2. Hence $h^0(Y,K_Y+\tB-2N_7)>0$ and the
restriction map $H^0(Y, K_Y+\tB-N_7)\to H^0(N_7,
\OO_{N_7}(2))$ is not surjective. This is a contradiction to Lemma
\ref{bir}, since $Z$ is an element of $\cN_0$.
\end{proof} Let $\cZ\subset (\pp^1\times\pp^1\times\pp^1)\times \cN$ be the universal
family. The group $G$ acts on
$\cZ$ preserving the fibres of
$\cZ\to\cN$, hence we can take the quotient and obtain a family
$\mathfrak{S}:=\cZ/G\to
\cN$, which is easily seen to be flat. We can also define a polarization $\cB$ on
$\kS$ as follows: we modify the $G-$action on the line bundle
$\OO_{\pp^1\times\pp^1\times\pp^1}(2,2,2)$ considered in Construction
\ref{construction} 
 by multiplying it with the nontrivial character of $G$ orthogonal to the subgroup
$G_0$. The effect of this choice is that the $G-$invariant sections now correspond to
the subspace
$T_1$. Denote by $Bs\subset \pp^1\times\pp^1\times\pp^1$ the base locus of the system
$|T_1|$. The restriction of $\OO_{\pp^1\times\pp^1\times\pp^1}(2,2,2)$ to  
$(\pp^1\times\pp^1\times\pp^1)\setminus Bs$ is generated at every point by global
sections which are invariant for the chosen linearization, hence it is the pull back
of a line bundle $\overline{\cB}$ from  the quotient
$((\pp^1\times\pp^1\times\pp^1)\setminus Bs)/G$. One can check that, by the
definition of $\cN$, the family $\cZ$ is contained in
$(\pp^1\times\pp^1\times\pp^1\setminus Bs)\times \cN$.  Hence the projection onto the
first factor induces  a map
$\kS\to (\pp^1\times\pp^1\times\pp^1\setminus Bs)/G$.  We let   $\cB$ be the pull
back of
$\overline{\cB}$ via this map.   For every
$t\in \cN$ the elements of
$T_1$ give global sections of
 the restriction
$B_t$ of
$\cB$ to the fibre $\Si_t$ of $\kS$ at $t$. By Construction \ref{construction} and 
Lemma
\ref{wiw},
$(\Si_t, B_t)$ is a pair as in the set-up of \S \ref{secex} and by Theorem
\ref{structure} all pairs
$(\Si,B)$ occur as $(\Si_t,B_t)$ for some $t\in N$.
\smallskip

A {\em simultaneous resolution} of a flat family $\cS \to U$ of normal projective
surfaces is a flat family $\cY\to U$ with a map $\cY\to \cS$ over $U$ such that for
every $t\in U$ the restricted  map $Y_t\to S_t$ is the minimal resolution of the
singularities of $S_t$.

\begin{prop}\label{resolution} The family $\kS\to \cN$ admits a simultaneous resolution $\cY\to
\cN$.
\end{prop}
\begin{proof} Consider the family $\cZ\subset ((\pp^1\times\pp^1\times\pp^1)\setminus
Bs)\times \cN$. Let $\hat{\pp}\to (\pp^1\times\pp^1\times\pp^1)\setminus Bs$ be the
blow up of the union of the fixed loci of the nonzero elements of $G$ and denote by
$\cZ'\subset 
\hat{\pp}\times \cN$ the pull back of $\cZ$. There is an induced $G-$action on $\cZ'$
with the property that the fixed loci of all the elements are (possibly empty) divisors. Denote by
$Disc\subset 
\hat{\pp}\times\cN$ the set of singular points of the fibres of $\cZ'$, i.e. the
inverse image of the incidence variety $\cI\subset 
(\pp^1\times\pp^1\times\pp^1)\times \cN$.  By the definition of $\cN$, the set  $\cI$ does
not meet the exceptional locus of the blow up, hence $Disc$ is isomorphic to $\cI$
and, in particular, it is smooth. The group 
$G$ acts freely on $Disc$ and the induced map $Disc\to\cN$ is the quotient by this
action.  By blowing up $Disc$ inside $ \hat{\pp}\times\cN$ and taking the strict transform pf $\cZ'$
one obtains a family
$\cZ''$ with an induced $G-$action.  The fibre of $\cZ''$ over a point $t\in \cN$ is
the blow up at the isolated fixed points of $G$ of the minimal resolution of the
fibre $Z_t$ of $\cZ$ at $t$. We set $\cY:=\cZ''/G$. Since  $G$ acts fibrewise, we
have an induced map $\cY\to\cN$. The family  $\cY\to \cN$ is smooth,  and hence flat,
since both the base $\cN$ and the fibres are smooth. By construction, there is a
natural map $\cY\to\kS$, commuting with the projections onto $\cN$,  which  
restricts to the minimal desingularization $Y_t\to\Si_t$ for every $t\in \cN$.
\end{proof}

\begin{lem}\label{section} The family of smooth surfaces $\cY\to\cN$ admits a section
$Sec\subset \cY$.
\end{lem}
\begin{proof} Consider the family $\cZ\subset (\pp^1\times \pp^1\times
\pp^1)\times\cN$. For every
$t\in
\cN$ the curve $C:=\{(1,1, 1,1)\}\times \pp^1$ meets $Z_t$ transversally at two
smooth points. Indeed, the intersection number of $C$ and $Z_t$ is equal to 2 and the
set $C\cap  Z_t$ is invariant under the action of the element $(e_2, e_2, e_2)$ of
$G$. Since $(e_2, e_2, e_2)$ acts freely on $Z_t$ by the definition of $\cN$,  the
only possibility is that $C\cap Z_t$ consists of two distinct points, that are not
fixed by any element of $G$. The intersection of
$\cZ$ with the subvariety
$C\times \cN $ is an
\'etale bisection of
$\cZ\to \cN$.  The image of this bisection in $\kS=\cZ/G$  is  a section of
$\kS\to\cN$ that intersects every fibre $\Si_t$ at a smooth point, and we take  its
inverse image in
$\cY$ as the required  section $Sec$.
\end{proof}

We recall briefly from \cite[Ch. 0, \S 5]{git}  (see also \cite{groth}) the main facts
about relative Picard schemes. 

Given a family $\cX\to T$, one  defines    the relative Picard functor   from the
category of schemes over $T$ to the category of sets. Given a scheme $T'\to  T$ , the
relative Picard functor associates to $T'$   the quotient of the group  of
isomorphism classes of line bundles on
$\cX\times _T T'$ by the subgroup of the classes of line bundles pulled back from
$T'$. If $\cX\to T$ admits a section, one can define the relative Picard functor also
by taking the isomorphism classes of ``normalized'' line bundles, namely of the line
bundles whose restriction to the pull back over $T'$ of the given section is trivial.
If $\cX\to T$ is flat and  projective with reduced irreducible fibres and it admits a
section, then the relative Picard functor is represented by a group scheme
$\Pic_{\cX/T}\to T$.  Therefore by Lemma \ref{section} we can consider the scheme
$\Pic_{\cY/\cN}\to \cN$, where $\cY$ is the family defined in Proposition \ref{resolution}.
Denote by
$\tcB$ the pull back to $\cY$ of the line bundle $\cB$ that we have previously defined
on $\kS$ and denote  by $Ex$ the exceptional divisor of the simultaneous resolution 
$\cY\to\kS$. The line bundle $\tcB\otimes\OO_{\cY}(Ex)$ defines a section $b\colon
\cN\to\Pic_{\cY/\cN}$. We define
$\tcN\subset \Pic_{\cY/\cN}$ to be the inverse image of $b(\cN)$ via the
multiplication by 2 map
$\Pic_{\cY/\cN}\to \Pic_{\cY/\cN}$. So $\tcN$ is closed in $\Pic_{\cY/\cN}$ and the
natural map $\tcN\to\cN$ is an \'etale double cover.  A point of $\tcN$ determines a
pair $(\Si, B)$ together with a solution $L\in\Pic(Y)$ of the linear equivalence 
$2L\equiv \tB+N_1+\dots+N_7$.

We are  finally ready to construct   varieties  that parametrize the isomorphism classes
of pairs
$(\Si,B)$ and the isomorphism classes of triples $(\Si, B,L)$ as above. We do this by
taking the  quotient of
$\cN$ and $\widetilde{\cN}$ by a suitable finite group. 

Let
$St(G)$ be the subgroup of $\Aut(\pp^1\times\pp^1\times\pp^1)$ consisting of the
elements 
$\ga$ such that
$\ga G\ga\inv=G$. It is easy to verify that $St(G)$ is a finite group. The group
$St(G)$  permutes  the $G-$eigenspaces $T_0,\dots, T_7$   of $H^0(\OO_{\pp^1\times
\pp^1\times\pp^1}(2,2,2))$. Since $T_0$ is the only eigenspace of dimension 5 and
$T_1$ is the only eigenspace of dimension 4, it follows that $St(G)$ preserves $T_0$
and $T_1$. In view of this observation, it follows from the definitions given so far
that $St(G)$ acts on
$\cN$, on the families
$\cZ$, $\kS$ and $\cY$ and  on the line bundle $\cB$ on $\kS$, and that all these
actions are compatible.   Clearly, the action of $St(G)$ on $\cY$  maps to itself 
the exceptional divisor
$Ex$ and therefore we also have an action of $St(G)$ on $\tcN$.

\begin{thm}\label{parametre}
\begin{enumerate}
\item The set of isomorphism classes of pairs $(\Si,B)$ as in the set up of\,  \S
\ref{secex} is in one-to-one correspondence with the quasiprojective variety
$\cN/St(G)$;
\item The set of isomorphism classes of triples $(\Si,B, L)$, where $(\Si, B)$ is a
pair as above and $L\in \Pic(Y)$ satisfies $2L\equiv \tB+N_1+\dots +N_7$,  is in
one-to-one correspondence with the quasi-projective variety $\tcN/St(G)$.
\end{enumerate}
\end{thm}
\begin{rem} One can formulate a suitable moduli problem for pairs $(\Si, B)$ and for
triples
$(\Si, B,L)$ and  it is very likely  that the spaces $\cN/St(G)$ and
$\tcN/St(G)$ are the corresponding coarse moduli spaces. Since we are mainly
interested in the applications to surfaces of general type with $p_g=0$, we will not
 pursue this any further.
\end{rem}

\begin{proof}[Proof of Theorem \ref{parametre}] We only  give the proof of (i), the
proof of (ii) being very similar. As we have already remarked, every pair $(\Si, B)$
is isomorphic to $(\Si_t, B_t)$ for some
$t\in \cN$. In addition, it is clear from the construction that if $\ga$ is an
element of
$St(G)$ then for every $t\in \cN$  the pairs corresponding to $t$ and $\ga t$ are
isomorphic.
 On the other hand, assume that $t, t'\in \cN$ give isomorphic pairs
$(\Si, B)$ and
$(\Si', B')$.  The $G-$covers $Z\to \Si$ and $Z'\to\Si'$ are defined intrinsically,
hence the isomorphism lifts to an isomorphism $Z\to Z'$. We notice that, up to the
ordering,  the  pencils
$|F_i|$ on
$\Si$ and the pencils $|F_i'|$ on $\Si'$ are  determined by
$B$, respectively $B'$. Indeed, at least one of the systems $|B|$ or $|K_{\Si}+B|$
maps $\Si$ birationally onto a singular  sextic in 
$\pp^3$ and the double fibres of the elliptic pencils $|F_i|$ are the inverse images
of the 6 double lines of the sextic (see \cite[Ch.IV, \S 9]{codo}). Hence, up to a permutation of the
factors of
$\pp^1\times\pp^1\times\pp^1$ (which is an element of $St(G)$), we may assume that for
$i=1,2,3$  the isomorphism $\Si\to \Si'$  maps $|F_i|$ to
$|F_i'|$. It follows from the proof of Theorem \ref{structure} that the isomorphism
$Z\to Z'$ is induced by an automorphism of $\pp^1\times\pp^1\times\pp^1$ compatible
with the $G-$action, namely by an element of
$St(G)$. 
\end{proof}

Since $\cN$ is irreducible by Theorem \ref{Ngeom},   the variety $\cN/St(G)$ is also
irreducible. The variety $\tcN$, being an \'etale double cover of $\cN$,   either is
irreducible or it is the disjoint union of two components isomorphic to $\cN$. We
close this section by showing that, in any case, taking the quotient of $\tcN$
by $St(G)$ we get an  irreducible variety. 

\begin{prop}\label{Nirred} The variety $\tcN/St(G)$ is irreducible.
\end{prop}
\begin{proof} The variety $\tcN$, being an \'etale cover of $\cN$,  is smooth by
Lemma \ref{Nsmooth} and thus $\tcN/St(G)$ is normal. So, to prove that $\tcN/St(G)$
is irreducible it suffices to show that it is connected. We do this by showing that
there exist a point $t$ in $\cN$ and an automorphism $\ga\in St(G)$ such that $\ga
t=t$ but  $\ga$ exchanges the two points of $\tcN$ lying over $t$. This amounts to
finding a pair $(\Si,B)$ such that there exists an automorphism $h$ of $\Si$ with
$h^*B\equiv B$ and such that the induced automorphism of $Y$ exchanges the two
solutions in $\Pic(Y)$ of the relation $2L\equiv \tB+N_1+\dots +N_7$. Indeed by Theorem
\ref{structure} such a pair is isomorphic to
$(\Si_t, B_t)$ for some $t$. Moreover, 
$h$ induces an automorphism
$h'$ of $Z$, since $Z$ is defined intrinsically. As we have observed in the proof of
Theorem \ref{parametre}, the set of pencils $|F_1|$, $|F_2|$, $|F_3|$ is determined
uniquely by $B$, hence it is preserved by $h$. It follows that $h'$ is compatible
with the embedding $Z\to |C_1|\times |C_2|\times |C_3|$ (cf. proof of Theorem
\ref{structure}, Step 2). In other words,  if we identify
$|C_1|\times |C_2|\times |C_3|$ with $\pp^1\times\pp^1\times \pp^1$ as
in the proof of Theorem \ref{structure}, then $h'$ is induced by  an element $\ga$ of
$St(G)$.  The pair $(\Si, B)$ that we construct is a special instance of  Example
\ref{ex1} (cf. also
\cite{mp3}).

So we let
$\Z_2^2=\{1, e_1, e_2, e_3\}$ act on   a product
$D_1\times D_2$ of elliptic curves by $(x,y)\stackrel{e_1}{\mapsto}(-x,y+b)$,
$(x,y)\stackrel{e_2}{\mapsto}(x+a,-y)$, where $a\in D_1$ and $b\in D_2$ are nonzero
elements of order 2. The quotient surface $\overline{\Si}$ is an Enriques surface
with $8$ nodes and has two elliptic pencils $|F_1|$, $|F_2|$  such that $F_1F_2=4$,
induced by the projections of
$D_1\times D_2$ onto the two factors. One of the double fibres of $|F_1|$ occurs over
the image in
$\pp^1=D_1/\Z_2^2$ of the points
$0$ and
$a$ and the other one occurs over the image of the remaining  $2-$torsion points
$a_1$ and
$a_2$. The fibres over the image in $\pp^1$ of the fixed points of $x\mapsto -x+a$
map contain 4 nodes each and they give rise two fibres of type $I_0^*$ on the
resolution $Y$ of
$\overline{\Si}$.  Now we assume in addition that
$D_1$ admits an  automorphism
$\tau$ of order 4 fixing the origin
$0$. The fixed locus of $\tau$ consists of the origin and of another point of order
$2$. Hence we may  take   $a$ in the above construction to be a fixed point of
$\tau$. We observe that $\tau$ exchanges the points $a_1$ and $a_2$. Consider the
automorphism
$h_0\colon D_1\times D_2\to D_1\times D_2$ defined by $(x,y)\mapsto (\tau x+a_1, y)$.
The automorphism  $h_0$  commutes with the elements of $\Z_2^2$, hence it induces an
automorphism $\bar{h}$ of the quotient surface $\overline{\Si}$, that clearly maps
each fibre of
$|F_2|$ to itself. The square of the map $x\to\tau x+a_1$ is equal to the map 
$x\mapsto -x+a$. Thus
$x\mapsto \tau x+a_1$ has order 4 and it fixes 2 points, that are necessarily also
fixed points of
$x\to-x+a$. Hence $\bar{h}$ maps to itself  each of the fibres with 4 nodes of
$|F_1|$ and it induces the identity on one of them.  On the other hand, $\bar{h}$
exchanges the two double fibres of
$|F_1|$. 
 We let
$\Si$ be the surface obtained by resolving one of the singular points of
$\overline{\Si}$ that are fixed by $\bar{h}$, we denote by $C$ the exceptional curve
of $\Si\to\overline{\Si}$ and we set
$B:=|F_1+F_2-C|$, where we omit to denote  pull backs. Clearly, $\bar{h}$ induces an
automorphism $h$ of $\Si$ and an automorphism of $Y$ that we  also denote by
$h$. As usual we denote by $|\tF_i|$, $i=1,2$,  the pull back of $|F_i|$ to $Y$ and by
$2\tE_i$, $2\tE_i'$ the double fibres of $|\tF_i|$. Furthermore we let 
$C_1$ and
$C_2$ be the  multiple components of the two  fibres of type
$I_0^*$ of $|\tF_1|$. Then the solutions in $\Pic(Y)$ of the relation $2L\equiv
\tB+N_1+\dots +N_7$ are the linear equivalence classes of $3\tE_1+\tE_2-C_1-C_2-C$
and of
$3\tE_1'+\tE_2-C_1-C_2-C$. It is clear by the above description that  these classes
are exchanged by $h$.
 \end{proof}
\section{A new family of  surfaces with
$p_g=0$ and $K^2=3$} \label{newfamily}

In this section we apply the previous results to the study of the moduli of surfaces
of general type with $p_g=0$ and $K^2=3$. We refer the reader to    M. Manetti's
Ph.D. thesis \cite{ma} for an excellent survey of the known results on this moduli
space.

We keep the notation from the previous sections.  Also  we let  
$\mathscr{M}$ be the moduli space of (canonical models) of surfaces of general type
with $p_g=0$ and $K^2=3$, and we denote by $\mcE$ the subset of
$\mcM$ consisting of the canonical surfaces whose bicanonical map  is composed
with an involution $\si$ such that the quotient surface $X/\si$ is birational to   an
Enriques surface. Notice that, if $X$ belongs to $\mcE$, then, by Theorem
\ref{gentoen},  $X/\si$ is in fact a nodal Enriques surface with 7 nodes.

\begin{thm}\label{irrmoduli}  The set $\mcE$ is constructible. \newline The closure
$\overline{\mcE}$ of
$\mcE$ in $\mcM$ is  irreducible and uniruled  of dimension $6$.
\end{thm}
\begin{proof} Let $\cN$, $\tcN$ be the spaces introduced in \S \ref{secmod} and let
$\tq\colon \tcY\to \tcN$ be the family obtained by pulling back the family
$\cY\to\cN$ defined in Proposition \ref{resolution}. We denote again by  $\tcB$ and $Ex$  the
pullbacks on
$\tcY$ of the corresponding objects of $\cY$. By Lemma \ref{section}, the family 
$\cY\to\cN$ has a section $Sec$, that induces a section of
$\tq$ that we denote again  by $Sec$.   Up to tensoring with a line bundle pulled
back from $\widetilde{\cN}$, we may also assume that the line bundle
$\tcB$ is normalized with respect to the section $Sec$, namely that its restriction
to $Sec$ is trivial. Then, if we denote by $\cL$ the pull back  to
$\tcY$ of the normalized Poincar\'e line bundle on $\cY\times_{\cN}\Pic_{\cY/\cN}$,
we have the equivalence relation $2\cL\equiv \tcB+Ex$. By the semicontinuity
theorem,  the sheaf
$\tq_*\tcB$ is a rank 4 vector bundle on $\widetilde{\cN}$.  We let $\cV$ be the
total space of
$\tq_*\tcB$ and we consider the family
$\widetilde{\cY}\times_{\widetilde{\cN}}\cV\to
\cV$.  We denote again by
 $\tcB$ and $Ex$ the pull backs to
$\widetilde{\cY}\times_{\widetilde{\cN}}\cV$ of the corresponding line
bundles/divisors of $\tcY$. Let $Taut$ be the zero locus  on
$\widetilde{\cY}\times_{\widetilde{\cN}}\cV$ of the tautological section of
$\tcB$ and let $q_1\colon \cV\to {\widetilde{\cN}}$ be the natural map.
 The fibre at
$v\in \cV$ of $Taut\to \cV$
 is a curve $B_v$ contained   the fibre $Y_{q_1(v)}$ of
$\widetilde{\cY}\to{\widetilde{\cN}}$ at $q_1(v)$. We say that $v$ is admissible if 
$B_v$ is contained in the smooth part of $Y_{q_1(v)}$ and it has at most negligible
singularities. We denote by
$\cV_{ad}\subset \cV$ the set of admissible points. It is easy to see that
$\cV_{ad}$ is a dense open subset of $\cV$. We denote by $\cS$ the double cover of
$\widetilde{\cY}\times_{\widetilde{\cN}}\cV_{ad}$ given by the relation $2\cL\equiv
Taut+Ex$. By Theorem \ref{entogen}, the fibre $S_v$ of $\cS\to\cV_{ad}$ over a point
$v\in \cV_{ad}$ is a  surface with canonical singularities, whose canonical model
belongs  to $\mcE$. Conversely, by Theorem \ref{gentoen}  every surface of $\mcE$ is
the canonical model of some surface $S_v$. Finally,  we denote by $p\colon \cX\to
\cV_{ad}$ the   relative canonical model of $\cS\to\cV_{ad}$. By a result of Iitaka
(\cite{Ii}) the plurigenera $p_m(S_v)$ are constant as functions of
$v$. Hence
$p_*\omega_{\cX/\cV_{ad}}^m$ is a vector bundle and, by the results of Bombieri, for
$m\ge 5$ the relative $m-$canonical map embeds $\cX$ into
$\pp(p_*\omega_{\cX/\cV_{ad}}^m)$. This means  that locally over $\cV_{ad}$ one  can
realize 
$\cX$ 
 as a family of subvarieties of a fixed projective space with constant Hilbert
polynomial. Since the space  $\cV_{ad}$ is smooth by Lemma \ref{Nsmooth}, the family 
$\cX\to\cV_{ad}$ is flat by
\cite[Ch. III, Thm. 9.9]{hartshorne}.
 By the properties of moduli spaces, the family  $\cX$ induces a morphism
$\Psi\colon\cV_{ad}\to\mcM$ whose image is
$\mcE$. This proves that $\mcE$ is constructible.  The map $\Psi$ factorizes through the natural map $\cV_{ad}\to
\pp(\cV)$. The image
$\overline{\cV_{ad}}$ of $\cV_{ad}$ in $\pp(\cV)$ is a dense open set. Arguing as in
the proof of Theorem
\ref{parametre} one sees that the action of
$St(G)$ on $\tcN$ lifts to an action on $\overline{\cV_{ad}}$ such that the quotient
$\overline{\cV_{ad}}/St(G)$ parametrizes the isomorphisms classes of pairs
$(X,\si)$.   Hence the map
$\overline{\cV_{ad}}\to \mcE$ factorizes through the quotient map
$\overline{\cV_{ad}}\to
\overline{\cV_{ad}}/St(G)$ and the latter map has finite fibres,  since a surface of
general type $X$ has finitely many automorphisms. The fibres of the  natural map
$\overline{\cV_{ad}}/St(G)\to
\tcN/St(G)$ are  irreducible of dimension 3, hence $\overline{\cV_{ad}}/St(G)$ is
irreducible of dimension 6 by Proposition \ref{Nirred} and it is uniruled, since
$\overline{\cV_{ad}}$ is an open set in a $\pp^3-$bundle. This remark completes the
proof.
\end{proof}
\begin{rem} With  some extra   work, one can show that, given a surface $X\in \mcE$,
the involution $\si$ is uniquely determined. Hence the map
$\overline{\cV_{ad}}/St(G)\to
\mcE$ in the proof of Theorem \ref{irrmoduli} is actually bijective.   
\end{rem}
\begin{rem} The main question left open by Theorem \ref{irrmoduli} is whether
$\overline{\mcE}$ is an irreducible component of $\mcM$. To answer this question one
has to consider for $X\in \mcE$ the natural map of functors $Def(X,\si)\to Def(X)$,
where $Def(X)$ denotes deformations of $X$ and $Def(X,\si)$ denotes deformations of
$X$ with an involution extending $\si$.  One needs to  decide whether this map is
surjective for a general $X$. To show that this is indeed the case,   it is enough to
exhibit one  surface
$X\in\mcE$ such that the map 
$Def(X,\si)\to Def(X)$ is smooth, and  smoothness can  in turn be checked by means
of  an infinitesimal computation. Unfortunately, although we can show  that for a
smooth  $X\in
\mcE$  the functor
$Def(X,\si)$ is smooth, we have not been able to  prove the  smoothness of
$Def(X,\si)\to Def(X)$. Notice that,  since the expected
dimension of
$Def(X)$ is equal to 4, Theorem \ref{irrmoduli} implies that the obstruction space  
$T^2_X$ of $Def(X)$ has dimension $\ge 2$ at every point of $\mcE$.
\end{rem}
\begin{cor} \label{fgroup} Let $S$ be a smooth surface such that the canonical model
of $S$ is in $\mcE$. Then:
$$\pi_1(S)\simeq \Z_2^2\times\Z_4.$$
\end{cor}
\begin{proof}Since blowing up does not change the fundamental group of a smooth
surface, we may assume that $S$ is minimal.
 By Theorem \ref{irrmoduli}, all the minimal surfaces whose canonical model is   in
$\mcE$ have the same fundamental group, so the statement follows by \cite[Thm.
3.1]{naie}.
\end{proof}

\begin{prop} \label{deg4} If $X\in \mcE$, then the bicanonical map $\fie$ of $X$ is a
morphism of degree either 2 or 4. The subset ${\mcE}_{d4}$ consisting of the surfaces
for which $\deg\fie=4$ is a closed subset of $\mcE$ of codimension 1 and its closure
$\overline{{\mcE}_{d4}}$ is irreducible.
\end{prop}
\begin{proof} The fact that  $\fie$ is a morphism of degree 2 or 4 is immediate by
Proposition
\ref{gentoen}. Since  $\deg\fie$ is a semicontinuous function of $X\in \mcM$, the
set
${\mcE}_{d4}$ is clearly closed in $\mcE$ and it is a proper subset of $\mcE$ by
Examples
\ref{ex2} and \ref{ex3}. To show the last part of the statement one proceeds as in
the proof of Theorem \ref{irrmoduli} by constructing a $5-$dimensional family of
surfaces that maps onto
${\mcE}_{d4}$ with finite fibres. By Theorem \ref{B+fibres}, Propositions
\ref{gentoen} and
\ref{entogen}, the fibres of this  family are  the  double covers $X\to \Si$, with $\Si$
an Enriques surface with 7 nodes, branched on the nodes and on a divisor $B$ with
negligible singularities, and such that the pair
$(\Si, B)$ is as in Example
\ref{ex1}. We omit the explicit construction of this family, which is standard by the
classification of Enriques surfaces with 8 nodes given in \cite{mp3}.
\end{proof}

\bigskip

\begin{tabbing} 1749-016 Lisboa, PORTUGALxxxxxxxxx\= 56127 Pisa, ITALY \kill
Margarida Mendes Lopes \> Rita Pardini\\ CMAF \> Dipartimento di Matematica\\
 Universidade de Lisboa \> Universit\a`a di Pisa \\ Av. Prof. Gama Pinto, 2
\> Via Buonarroti 2\\ 1649-003 Lisboa, PORTUGAL \> 56127 Pisa, ITALY\\
mmlopes@ptmat.ptmat.fc.ul.pt \> pardini@dm.unipi.it
\end{tabbing}


\begin{thebibliography}{BPV}
\bibitem[BPV]{bpv} W.~Barth, C.~Peters, A.~Van de Ven, {\em Compact complex
surfaces}, Ergebnisse der Mathematik und ihrer Grenzgebiete, {\bf3.} Folge, Band
{\bf4}, Springer-Verlag, Berlin (1984).
\bibitem [Be] {bv} A.~Beauville, {\em Complex Algebraic Surfaces}, Cambridge
University Press (1983).
\bibitem[Bo] {bo} E.~Bombieri, {\em Canonical models of surfaces of general type},
Publ. IHES {\bf 42}
 (1973), 447--495.
\bibitem [Bu] {bu} P.~Burniat,  {\em Sur les surfaces de genre
$P_{12}>0$},   Ann. Mat. Pura Appl., (4)  71  (1966).


\bibitem[BW]{bw} D. M. Burns, J.M.  Wahl, {\em Local contributions to global
deformations of surfaces}, Invent. Math. {\bf 26} (1974), 67--88.

\bibitem [Ca]{quatro} F.~Catanese,  {\em Singular bidouble covers and the
construction of interesting algebraic surfaces}, in {\em Algebraic Geometry:
Hirzebruch $70$},  A.M.S. Contemporary Mathematics, vol. 241 (1999), 97--120.




\bibitem[Ci] {ciro} C.~Ciliberto, {\it  The bicanonical map for surfaces of general
type}, Proceedings of Symposia in Pure Mathematics  {\bf 62.1}  (1997), 57--84.

\bibitem[Co]{cossec} F. Cossec, {\em Projective models of Enriques surfaces}, Math.
Ann. {\bf 265} (1983), 283-334.
\bibitem[CD]{codo} F. Cossec, I. Dolgachev, {\em Enriques surfaces I}, Progress in
Mathematics {\bf 76}, Birkh\"auser (1989).

\bibitem[DMP]{nodes} I.~Dolgachev, M.~Mendes Lopes, R.~Pardini, {\em Rational
surfaces with many nodes}, Compositio Math. {\bf 132} (2002), no. 3, 349--363.

\bibitem[Gr]{groth} A. Grothendieck, {\em Fondements de la G\'eom\'etrie
Alg\'ebrique}, expos\'es 232, 236, collected Bourbaki talks, Paris 1962.



\bibitem[Ha]{hartshorne} R.~Hartshorne, {\em Algebraic Geometry}, G.T.M. 52,
Springer-Verlag, New York (1977).

\bibitem[In]{inoue} M.~Inoue, {\em Some new surfaces of general type}, Tokyo J. of
Math., {\bf 17}  (2) (1994), 295--319.

\bibitem[Ii]{Ii} S. Iitaka, {\em Deformations of compact complex surfaces II}, J.
Math. Soc. Japan {\bf 22} (1970), 247--261.

\bibitem[Ke]{keum} J.H.Keum, {\em Some new surfaces of general type with $p_g=0$},
(preprint 1988). 

 \bibitem[Li]{lint} J.H. van Lint, {\em Introduction to coding theory}, Second
edition,  Graduate Texts in Mathematics {\bf 86}, Springer-Verlag, Berlin, 1992.

\bibitem[Ma]{ma} M. Manetti, {\em Degeneration of algebraic surfaces and applications
to moduli problems}, Tesi di Perfezionamento, Scuola Normale Superiore -- Pisa (1996).
\bibitem[MP1]{mp} M.~Mendes Lopes, R.~Pardini, {\em The bicanonical map of surfaces
with
$p_g=0$ and $K^2\ge 7$},  Bull. London Math. Soc. {\bf 33} (2001), 1--10.

\bibitem[MP2]{mp2} M. Mendes Lopes, R. Pardini, {\em A connected component of the
moduli space of surfaces of general type with $p_g=0$}, Topology {\bf 40} (5)
 (2001),  977--991.

\bibitem[MP3]{mp3} M.~Mendes Lopes, R.~Pardini, {\em Enriques surfaces with eight
nodes}, Math. Z. {\bf 241} (2002), 673--683.
\bibitem[MP4]{bica2} M.~Mendes Lopes, R.~Pardini, {\em The bicanonical map of surfaces
with
$p_g=0$ and $K^2\ge 7$, II},  Bull. London Math. Soc. (to appear).
\bibitem[Mu]{git} D. Mumford, {\em Geometric invariant theory},  Ergebnisse der
Mathematik und
 ihrer Grenzgebiete, Neue Folge, Band {\bf30}, Springer-Verlag, Berlin (1965).
\bibitem[Na]{naie} D. Naie, {\em Surfaces d'Enriques et une construction de surfaces
de type g\'en\'eral avec $p_g=0$}, Math. Z. {\bf 215} (2) (1994), 269--280.
\bibitem[Pa1]{ritaabel} R. Pardini, {\em Abelian covers of algebraic varieties},   J.
reine angew. Math.  {\bf 417} (1991), 191--213.
\bibitem[Pa2]{pianidoppi} R. Pardini, {\em The classification of double planes of
general type with $p_g=0$ and $K^2=8$}, J. of Algebra {\bf 259} (2003), 95-118.
\bibitem[Pe]{peters} C.~Peters, {\em On certain examples of surfaces with $p_g=0$ due
to Burniat},  Nagoya Math. J., {\bf 166} (1977), 109--119.

\bibitem [Re] {red} I.~Reider, {\em Vector bundles of rank 2 and linear systems on
algebraic surfaces}, Ann. of Math., {\bf127} (1988), 309--316.

\bibitem [X1]{xiaocan} G.~Xiao, {\em Finitude de l'application bicanonique des
surfaces de type g\'en\'eral}, Bull. Soc. Math. France, {\bf113} (1985), 23--51.
\bibitem [X2]{xiao2} G.~Xiao, {\em Degree of the bicanonical map of a surface of
general type},  Amer. J. of Math. {\bf 112} (5) (1990), 713--737.

\end{thebibliography}
\end{document}